\documentclass[12pt,sort&compress]{article}
\usepackage{amsfonts}
\usepackage{bm}
\usepackage{makecell}
\usepackage{amssymb,amsmath}
\usepackage{amscd}
\usepackage{lipsum}
\usepackage{latexsym}
\usepackage{CJK}
\usepackage{longtable}
\usepackage{enumerate}
\def\g{\gamma}

\def\Vir{\hbox{Vir}}
\def\cl{\centerline}

\def\vs{\vspace*}

\def\Z{\mathbb{Z}}
\def\CC{\mathbb{C}}
\def\C{\mathbb{C}}

\def\QED{\hfill$\Box$}

\def\t{\tilde}

\def\pa{\partial}
\def\la{\lambda}
\def\HV{\mathcal{HV}}
\def\SV{\widetilde{\mathcal{SV}}}

\numberwithin{equation}{section}
\newtheorem{theo}{Theorem}[section]
\newtheorem{defi}[theo]{Definition}
\newtheorem{coro}[theo]{Corollary}
\newtheorem{lemm}[theo]{Lemma}

\newtheorem{rema}[theo]{Remark}

\newtheorem{example}[theo]{Example}
\newtheorem{proposition}[theo]{Proposition}
\makeatletter
\newcommand\blfootnote[1]{%
\begingroup
\renewcommand\thefootnote{}\footnote{#1}%
\addtocounter{footnote}{-1}%
\endgroup
}

\def\@biblabel#1{#1.~}

\makeatother

\begin{document}
\vs{10pt} \cl{\large {\bf Cohomology of the extended Schr\"odinger-Virasoro conformal algebra}}
\cl{Lipeng Luo$^{1}$, Henan Wu$^{2}$}
\cl{\small{$^{1}$School of Mathematical Sciences, Tongji University, Shanghai, 200092, China}}
\cl{\small{$^{2}$School of Mathematical Sciences, Shanxi University, Taiyuan, 030006, China}}
\cl{\small E-mail: luolipeng@tongji.edu.cn, wuhenan@sxu.edu.cn
 }\vs{6pt}
\small\parskip .005 truein \baselineskip 3pt \lineskip 3pt

\noindent{\bf Abstract:} All the basic cohomology groups and reduced cohomology groups of the extended Schr\"odinger-Virasoro conformal algebra with trivial coefficients are completely determined. In particular, we introduce the notion of
the relative cohomology of Lie conformal algebra, and then we can express our main results in a more concise manner.
\blfootnote{Corresponding author: Henan Wu (wuhenan@sxu.edu.cn).}

\vs{5pt}

\vs{5pt}


%\vs{18pt}
\section{Introduction}
%%%%%%%%%%%%%%%
The cohomology theory of Lie algebras is an important subject of research. The study of cohomology theory of Lie algebras was initiated by Chevalley and Eilenberg in \cite{CE}, which aims to compute the real cohomology of the underlying topological space of a compact connected Lie group in terms of its associated Lie algebra. Cohomology of infinite-dimensional Lie algebras is related to invariant differential operators, combinatorial identities, integrable systems, Riemannian foliations and cobordism theory in \cite{F}. The cohomology of the Jacobson-Witt Algebras with trivial coefficients were investigated in \cite{ShuYao}. Cohomology of Lie algebras of polynomial vector fields on the line over fields of characteristic $2$ were determined in \cite{WFV}. There are many applications in the cohomology theory of Lie algebras in both structure and representation theory of Lie algebras \cite{CE,DAS,JS,YNI}.

The notion of Lie conformal algebra is introduced by Kac in \cite{K1}. It is well known that Lie conformal algebra has a close relationship with Lie algebras, especially the Lie algebras satisfying some locality. And its structure theory and representation theory were comprehensively studied in \cite{DK,CK}.
In this paper, we consider the cohomology theory of Lie conformal algebras. The general cohomology theory of conformal algebras with coefficients in an arbitrary conformal module was systematically
developed in \cite{BKV}.
All cohomology groups of the Virasoro conformal algebra, the Heisenberg-Virasoro conformal algebra and the Schr\"odinger-Virasoro conformal algebra with trivial coefficients were determined in \cite{BKV,YW,WL}, respectively. Also, some low dimensional cohomology groups of the general Lie
conformal algebras $gc_N$ were determined in \cite{S}. There exists a correspondence between the cohomology of Lie conformal algebra and that of its annihilation Lie algebra \cite{BKV}.
So the study of cohomology theory of Lie conformal algebra helps us to understand the cohomology theory of Lie algebra, which seems to be more complicated.

In the present paper, we focus on the cohomology of the extended Schr\"odinger-Virasoro conformal algebra, denoted by $\SV$ and introduced in \cite{SY}. By definition, $\SV$ is a finite Lie conformal algebra with free generators $L,N,Y,M$ subject to the
following nontrivial $\lambda$-brackets\vspace*{-7pt}
\begin{eqnarray}
&&[L_\lambda L]=(\partial+2\lambda)L,\ \ \ [Y_\la L]=(\frac12\pa+\frac32\la)Y, \label{lamda-bracket11}\\
&&{[L_\lambda Y]}=(\partial+\frac32\lambda)Y,\ \ \,{[Y_\lambda Y]}=(\partial+2\lambda)M,\label{lamda-bracket22}\\[3pt]
&&{[L_\lambda M]}=(\partial+\lambda)M,\ \ \,[M_\la L]=\la M,\label{lamda-bracket33} \\[5pt]
&&[L_\lambda N]=(\partial+\lambda)N,\ \ \ \,[N_\lambda Y]=Y, \ \ \,[N_\lambda M]=2M.
\end{eqnarray}
The extended Schr\"odinger-Virasoro conformal algebra $\SV$ contains many interesting examples as its subalgebras. For instance,
$\SV$ contains the
Virasoro conformal algebra $Vir=\mathbb{C}[\partial]L$,
the
Heisenberg-Virasoro conformal algebra $\mathcal{HV}=\mathbb{C}[\partial]L\bigoplus
\mathbb{C}[\partial]M$ and the
Schr\"odinger-Virasoro conformal algebra $\mathcal{SV}=\mathbb{C}[\partial]L\bigoplus
\mathbb{C}[\partial]Y\bigoplus
\mathbb{C}[\partial]M$, whose whose representation theory and cohomology theory were investigated in \cite{CK,BKV,WY,YW,WL}, respectively.

%%%%%%%%%%%%%
The first relative cohomology group of the Lie algebra of smooth vector fields with coefficients in the space of trilinear differential operators that acts on tensor densities were determined in \cite{BAL}. Evens and Graham studied the Belkale-Kumar family of cup products on the cohomology of a generalized flag variety and gave an alternative construction of the family using relative Lie algebra cohomology(see\cite{EG}). Muzere proposed a vanishing property for the Hochschild relative cohomology groups, which was revealed that the relative cohomology groups $H^*(g,f,V)$ and $H^*(Ug, Uf, V)$ are not in general isomorphic(see \cite{MM}). It is worth mentioning that we introduce the notion of relative cohomology of Lie conformal algebra in the present paper, which can be viewed as a generalization of the relative cohomology of Lie algebra \cite{F}.
We believe that this notion can help us to understand the cohomology of the subalgebra and the quotient algebra of a given Lie conformal algebra.
This is a motivation of presenting our work.

The rest of the paper is organized as follows. In Section 2, we recall some basic definitions, notations, and related known results about Lie conformal algebras. In particular, we introduce the notion of relative cohomology of Lie conformal algebra. %Then the main theorem of this paper is given.
In Section 3, we determine the basic cohomology groups of the extended Schr\"odinger-Virasoro conformal algebra with coefficients in its trivial module $\mathbb{C}_a$. In Section 4, we compute the reduced cohomology groups of the extended Schr\"odinger-Virasoro conformal algebra with coefficients in its module $\mathbb{C}_a$. In Section 5, we investigate the basic relative cohomology and reduced relative cohomology of the extended Schr\"odinger-Virasoro conformal algebra $\SV$ modulo $\mathcal{B}$ with coefficients in a trivial module $\mathbb{C}_a$, where $\mathcal{B}$ is a subalgebra of $\SV$.

Our main results are summarized in Theorems \ref{th1}, \ref{th2}, \ref{th3}.%, 5.4, 5.5, 5.6, 5.8, 5.11.

Throughout this paper, we use notations $\mathbb{C}$, $\mathbb{Z}$ and $\mathbb{Z_{+}}$ to represent the set of complex numbers, integers and nonnegative integers, respectively. In addition, all vector spaces and tensor products are over $\mathbb{C}$. In the absence of ambiguity, we abbreviate $\otimes_{\mathbb{C}}$ to $\otimes$.
%%%%%%%%%%%%%%%
\vs{8pt}
\section{Preliminaries}

   In this section, we recall some basic definitions, notations and related results about Lie conformal algebras for later use. For a detailed description, one can refer to \cite{BKV,CK,DK,YW}. %Then we present our main results.
\subsection{Lie conformal algebra}
\begin{defi}
\begin{em}(\cite{DK})
A Lie conformal algebra $\mathcal {A}$ is a $\CC[\partial ]$-module endowed with a $\CC$-bilinear map $\mathcal {A}\otimes \mathcal {A}\rightarrow \CC[\lambda]\otimes \mathcal {A}$, $a\otimes b \mapsto [a_\lambda b]
$ subject to the following relations ($a, b, c\in \mathcal {A}$):
\begin{align*}
[\partial a_\lambda b]&=-\lambda[a_\lambda b],\ \ [ a_\lambda \partial b]=(\partial+\lambda)[a_\lambda b] \ \ \mbox{(conformal\  sesquilinearity)},
\\
{[a_\lambda b]} &= -[b_{-\lambda-\partial}a] \ \ \mbox{(skew-symmetry)},
\\
{[a_\lambda[b_\mu c]]}&=[[a_\lambda b]_{\lambda+\mu
}c]+[b_\mu[a_\lambda c]]\ \ \mbox{(Jacobi \ identity)}.
\end{align*}
A conformal algebra is called \emph{finite} if it is finitely generated as a $\C[\partial]$-module.
\end{em}
\end{defi}

\subsection{Conformal module}
\begin{defi}
\begin{em}(\cite{CK})
A conformal module $V$ over a Lie conformal algebra $\mathcal {A}$
is a $\mathbb{C}[\partial]$-module equipped with a $\CC$-bilinear map
$\mathcal {A}\otimes V\rightarrow
V[\lambda]$, $a\otimes v\mapsto a_\lambda v$, satisfying the following relations for any $a,b\in\mathcal {A}$, $v\in V$,
\begin{eqnarray*}
&&a_\lambda(b_\mu v)-b_\mu(a_\lambda v)=[a_\lambda b]_{\lambda+\mu}v,\\
&&(\partial a)_\lambda v=-\lambda a_\lambda v,\ a_\lambda(\partial
v)=(\partial+\lambda)a_\lambda v.
\end{eqnarray*}
If $V$ is finitely generated over $\mathbb{C}[\partial]$, then $V$ is simply called \emph {finite}.
\end{em}
\end{defi}

\begin{example} Let $\mathcal{A}$ be an arbitrary Lie conformal algebra and $a\in\C$.  Then $\mathcal{A}$ admits a family of $1$-dimensional modules $\CC_a$ defined by
$$\CC_a=\CC,\ \partial v=a v,\ \mathcal {A}_\lambda v=0,\ \forall v\in\CC_a.$$
And we abbreviate $\C_0$ to $\C$ in the sequel. It is easy to check that the modules $\mathbb{C}_a$ with $a \in \mathbb{C}$ exhaust all trivial irreducible $\mathcal{A}$-modules.

\end{example}
%%%%%%%%%%%%

\subsection{Basic cohomology}
\begin{defi}\label{cochain}\rm (\cite{BKV}) An \emph{$n$-cochain} ($n\in\Z_+$) of a Lie conformal algebra $\mathcal{A}$ with coefficients in an
$\mathcal{A}$-module $V$ is a $\CC$-linear map\vs{-5pt}
\begin{eqnarray*}
\gamma:\mathcal{A}^{\otimes n}\rightarrow V[\la_1,\cdots,\la_n],\ \
\ a_1\otimes\cdots \otimes a_n \mapsto
\g_{\la_1,\cdots,\la_n}(a_1,\cdots,a_n)
\end{eqnarray*}
satisfying the following conditions:\begin{itemize}\parskip-3pt
\item[\rm(1)] $\g_{\la_1,\cdots,\la_n}(a_1,\cdots,\pa a_i,\cdots,
a_n)=-\la_i\g_{\la_1,\cdots,\la_n}(a_1,\cdots, a_n)$ \ (conformal antilinearity),
\item[\rm (2)] $\g$ is skew-symmetric with respect to simultaneous permutations
of $a_i$'s and $\la_i$'s \ (skew-symmetry).
\end{itemize}
\end{defi}

%As usual, let $\mathcal{A}^{\otimes 0}= \CC$, so that a $0$-cochain
%is an element of $V$.
Denote by  ${\t C}^n(\mathcal {A},V)$ the set
of all $n$-cochains. The differential $d_n$ of an $n$-cochain $\g$ is
defined as follows:
\begin{align}\label{ddd}
&(d_n\g)_{\la_1,\cdots,\la_{n+1}}(a_1,\cdots,a_{n+1})\nonumber\\
&=\mbox{$\sum\limits_{i=1}^{n+1}$}(-1)^{i+1}a_{i_{\la_i}}\g_{\la_1,\cdots,\hat{\la_i},\cdots,\la_{n+1}}(a_1,\cdots,\hat{a_i},\cdots,a_{n+1})\nonumber\\
&+\mbox{$\sum\limits_{1\le i<j\le n+1}$}(-1)^{i+j}\g_{\la_i+\la_j,\la_1,\cdots,\hat{\la_i},\cdots,\hat{\la_j},\cdots,\la_{n+1}}([a_{i_{\la_i}}a_j],a_1,\cdots,\hat{a_i},\cdots,\hat{a_j},\cdots,a_{n+1}),
\end{align}
where $\g$ is linearly extended over the polynomials in $\la_i$. %In
%particular, if $\g\in V$ is a $0$-cochain, then
%$(d_0\g)_\la(a)=a_\la\g$.

It was shown in \cite{BKV} that the operator $d$ preserves the space of cochains and $d^2=0$. Thus the cochains of a Lie conformal algebra $\mathcal{A}$ with coefficients in an $\mathcal{A}$-module $V$ form a complex, called the {\it basic complex} :
\begin{equation}
\begin{CD}
\cdot\cdot\cdot@>>> \t C^{n-1}(\mathcal{A},V)  @>{\rm d}_{n-1}>> \t C^{n}(\mathcal{A},V) @>{\rm d}_n>>\t C^{n+1}(\mathcal{A},V)  @>>>\cdot\cdot\cdot\ .
\end{CD}
\end{equation}
%which will be denoted by ${\t C}^\bullet(\mathcal {A},V)=\bigoplus_{n\in \mathbb{Z}_+}{\t C}^n(\mathcal {A},V)$.
The related cohomology is called the {\it basic cohomology} of the Lie conformal algebra $\mathcal{A}$ with coefficients
in its module $V$ and denoted by ${\rm \t H}^q (\mathcal{A},V),\ q\in\Z_+$.%, for more details see the following definition.

\begin{defi}
\begin{em}
An element $\gamma$ in ${\tilde C}^q(\mathcal{A},V)$ is called
a {\it $q$-cocycle} if $d(\gamma)=0$; a {\it $q$-coboundary} if there exists a $(q-1)$-cochain $\phi\in\tilde
C^{q-1}(\mathcal{A},V)$ such that $\gamma=d(\phi)$. Two cochains
$\gamma_1$ and $\gamma_2$ are called {\it equivalent} if $\gamma_1-\gamma_2$ is a coboundary. \end{em}
\end{defi}

Denote by $\tilde D^q(\mathcal{A},V)$ and $\tilde B^q(\mathcal{A},V)$ the spaces
of $q$-cocycles and $q$-boundaries, respectively. Then, we can obtain that
\begin{eqnarray*}
{\rm \tilde H}^q(\mathcal{A},V)=\tilde D^q(\mathcal{A},V)/\tilde B^q(\mathcal{A},V)=\{\mbox{equivalent classes of
$q$-cocycles}\}.
\end{eqnarray*}

\subsection{Reduced cohomology}
One can define a (left) $\CC[\partial]$-module structure on $\t
C^n(\mathcal{A},V)$ by
\begin{eqnarray*}
(\partial\g)_{\la_1,\cdots,\la_n}(a_1,\cdots,
a_n)=(\partial_V+\mbox{$\sum\limits_{i=1}^n$}\la_i)\g_{\la_1,\cdots,\la_n}(a_1,\cdots,
a_n),
\end{eqnarray*}
where $\partial_V$ denotes the action of $\partial$ on $V$. Then $d\partial=\partial d$
and %thus $\partial{\t C}^\bullet(\mathcal {A},V)\subset {\t C}^\bullet(\mathcal {A},V)$, i.e.,
\begin{equation}
\begin{CD}
\cdot\cdot\cdot@>>>\partial  \t C^{n-1}(\mathcal{A},V)  @>{\rm d}_{n-1}>>\partial  \t C^{n}(\mathcal{A},V) @>{\rm d}_n>>\partial \t C^{n+1}(\mathcal{A},V)  @>>>\cdot\cdot\cdot
\end{CD}
\end{equation}
forms a subcomplex of the basic complex. The quotient complex %$C^\bullet(\mathcal {A},V)={\t C}^\bullet(\mathcal {A},V)/\partial{\t C}^\bullet(\mathcal {A},V)=\bigoplus_{n\in \mathbb{Z}_+}C^n(\mathcal {A},V)$, i.e.,
\begin{equation}
\begin{CD}
\cdot\cdot\cdot@>>>
\frac{\t C^{n-1}(\mathcal{A},V)}{\partial \t C^{n-1}(\mathcal{A},V)}
@>{\bar{\rm d}}_{n-1}>>
\frac{\t C^{n}(\mathcal{A},V)}{\partial \t C^{n}(\mathcal{A},V)}
@>{\bar{\rm d}}_n>>
\frac{\t C^{n+1}(\mathcal{A},V)}{\partial \t C^{n+1}(\mathcal{A},V)}  @>>>\cdot\cdot\cdot
\end{CD}
\end{equation}
%\begin{eqnarray*}
%\rightarrow\frac{\t C^{n-1}(\mathcal{A},V)}{\partial \t C^{n-1}(\mathcal{A},V)}\rightarrow \frac{\t C^{n}(\mathcal{A},V)}{\partial \t C^{n}(\mathcal{A},V)}\rightarrow\frac{\t C^{n+1}(\mathcal{A},V)}{\partial \t C^{n+1}(\mathcal{A},V)}\rightarrow
%\end{eqnarray*}
is called the {\it reduced complex}. And its cohomology is called the \emph{reduced
cohomology} of the Lie conformal algebra $\mathcal{A}$ with coefficients
in $V$ and denoted by ${\rm H}^q (\mathcal{A},V),\ q\in\Z_+$.

\begin{rema}
The basic cohomology ${\rm \t H}^q (\mathcal{A},V)$ is naturally a $\CC[\partial]$-module, whereas the reduced cohomology ${\rm H}^q (\mathcal{A},V)$ is a complex vector space.
\end{rema}
The exact sequence $0\longrightarrow \partial\t C^\bullet \stackrel{i}{\longrightarrow} \t C^\bullet \stackrel{p}{\longrightarrow}C^\bullet \longrightarrow 0$ gives a long exact sequence of the cohomology groups:
\begin{align}\label{longexact}
0\longrightarrow& H^0(\partial\t C^\bullet) \stackrel{i_0}{\longrightarrow} \t H^0(\mathcal{A},V) \stackrel{p_0}{\longrightarrow} H^0 (\mathcal{A},V) \longrightarrow\nonumber\\
\longrightarrow& H^1(\partial\t C^\bullet) \stackrel{i_1}{\longrightarrow} \t H^1(\mathcal{A},V) \stackrel{p_1}{\longrightarrow} H^1 (\mathcal{A},V) \longrightarrow\\
\longrightarrow& H^2(\partial\t C^\bullet) \stackrel{i_2}{\longrightarrow} \t H^2(\mathcal{A},V) \stackrel{p_2}{\longrightarrow} H^2 (\mathcal{A},V) \longrightarrow\cdots.\nonumber
\end{align}

%Inspired by \cite{CK}, we can obtain the following result, which plays an important role in the classification of the reduced cohomology of  Lie conformal algebra $\mathcal{A}$ with coefficients in $V$.
\begin{proposition}(\cite{BKV})\label{pro2.12}
In degrees $\geq1$, the complexes $\t C^\bullet$ and $\partial\t C^\bullet$ are isomorphic under the map
\begin{align}
\t C^\bullet \to \partial\t C^\bullet, \quad \g \mapsto \partial\cdot \g.
\end{align}
Therefore, $H^q(\partial\t C^\bullet)\cong\t H^q(\mathcal{A},V)$ for $q\geq1$.
\end{proposition}

\begin{rema}
%The above proposition does not imply that
In the long exact sequence (\ref{longexact}),
the maps $H^q(\partial\t C^\bullet)\to\t H^q(\mathcal{A},V)$ induced by the embedding $\partial\t C^\bullet\subset\t C^\bullet$ are not isomorphisms.
\end{rema}

\subsection{Relative cohomology} Inspired by the notion of the relative cohomology of Lie algebra \cite{F}, we introduce the analogous notion of the relative cohomology of Lie conformal algebra.
%In order to study the properties of cohomology of a Lie conformal algebra $\mathcal{A}$ with coefficients in an
%$\mathcal{A}$-module $V$ from another angle, we can introduce the following definitions.

\begin{defi}\label{rcochain}\rm
Let $\mathcal{A}$ be a Lie conformal algebra and $\mathcal{B}$ a subalgebra of $\mathcal{A}$. Suppose that $V$ is an $\mathcal{A}$-module. Denote by ${\tilde C}^q(\mathcal{A};\mathcal{B},V)$ the subspace of the space ${\tilde C}^q(\mathcal{A},V)$, consisting of $q$-cochains $\g$ such that $\g_{\la_1,\cdots,\la_q}(a_1,\cdots, a_q)=(d\g)_{\la_1,\cdots,\la_{q+1}}(b_1,\cdots, b_{q+1})=0$ for any $a_1,b_1\in\mathcal{B}$.
Elements of the space ${\tilde C}^q(\mathcal{A};\mathcal{B},V)$ are called \emph{relative cochains} of $\mathcal{A}$ modulo $\mathcal{B}$ with coefficients in an
$\mathcal{A}$-module $V$.
\end{defi}

\begin{rema} Assume $U$ and $V$ are conformal modules over a Lie conformal algebra $\mathcal{A}$, then $U\bigotimes V$ becomes an $\mathcal{A}$-module under the following action
$$a_\lambda(u\otimes v)=(a_\lambda u)\otimes v+u\otimes (a_\lambda v),\ \ \partial(u\otimes v)=(\partial u)\otimes v+u\otimes (\partial v).$$
This definition is different from the notion of "tensor product" introduced in \cite{L}. But it coincides with the usual "tensor product" of modules of Lie algebra. %$U\bigotimes V$ is no longer finite.
So the usual n-th tensor power $V^{n}:=V\bigotimes\cdots \bigotimes V$ and the usual n-th exterior power $\Lambda^n V$ are also $\mathcal{A}$-modules.
If $\mathcal{B}$ is a subalgebra of $\mathcal{A}$, then $\Lambda^q(\mathcal{A}/\mathcal{B})$ is a module over $\mathcal{B}$.
And one can check that ${\tilde C}^q(\mathcal{A};\mathcal{B},V)=Hom_{\mathcal{B}}(\Lambda^q(\mathcal{A}/\mathcal{B}),V)$, which can be viewed as the equivalent definition of Definition \ref{rcochain}.
\end{rema}

Obviously, $d{\tilde C}^q(\mathcal{A};\mathcal{B},V)\subset{\tilde C}^{q+1}(\mathcal{A};\mathcal{B},V)$, so that relative cochains constitute a subcomplex of the complex ${\tilde C}^\bullet(\mathcal{A},V)$. This subcomplex is denoted by ${\tilde C}^\bullet(\mathcal{A};\mathcal{B},V)$, and its cohomology is called a \emph{basic relative cohomology} of the Lie conformal algebra $\mathcal{A}$ modulo $\mathcal{B}$ with coefficients in an
$\mathcal{A}$-module $V$, and denoted by ${\rm \t H}^q (\mathcal{A};\mathcal{B},V)$. And the inclusion ${\tilde C}^q(\mathcal{A};\mathcal{B},V)\subset {\tilde C}^q(\mathcal{A},V)$ induces a homomorphism ${\rm \t H}^q (\mathcal{A};\mathcal{B},V)\rightarrow {\rm \t H}^q (\mathcal{A},V)$, which is not necessarily injective.

Since $\partial\g=(\partial_V+\mbox{$\sum\limits_{i=1}^n$}\la_i)\g\in {\tilde C}^\bullet(\mathcal{A};\mathcal{B},V)$ for any $\g\in {\tilde C}^\bullet(\mathcal{A};\mathcal{B},V)$, the {\it reduced relative cohomology} of the Lie conformal algebra $\mathcal{A}$ modulo $\mathcal{B}$ with coefficients in an
$\mathcal{A}$-module $V$, denoted by ${\rm H}^q (\mathcal{A};\mathcal{B},V)$, is also well defined. %Similarly, there exists a homomorphism ${\rm  H}^q (\mathcal{A};\mathcal{B},V)\rightarrow {\rm H}^q (\mathcal{A},V)$.

\begin{rema}\label{rela}
Note that if $\mathcal{B}$ is an ideal of $\mathcal{A}$, then $\Lambda^q(\mathcal{A}/\mathcal{B})$ is a trivial $\mathcal{B}$-module and $${\tilde C}^q(\mathcal{A};\mathcal{B},V)=Hom_{\mathcal{B}}(\Lambda^q(\mathcal{A}/\mathcal{B}),V)\cong Hom_{\C}(\Lambda^q(\mathcal{A}/\mathcal{B}),V^{\mathcal{B}})={\tilde C}^q(\mathcal{A}/\mathcal{B},V^{\mathcal{B}})$$ where $V^{\mathcal{B}}$ is the $\mathcal{B}$-invariants, i.e.,
\begin{align*}
V^{\mathcal{B}}=\{v\in V\  |\  b\cdot v\ =\ 0 \ for\ all\ b\in\mathcal{B}\}.
\end{align*}
In particular, in this case the differentials in the complexes ${\tilde C}^\bullet(\mathcal{A};\mathcal{B},V)$ and ${\tilde C}^\bullet(\mathcal{A}/\mathcal{B},V^{\mathcal{B}})$ also coincide so that ${\tilde H}^q(\mathcal{A};\mathcal{B},V)={\tilde H}^q(\mathcal{A}/\mathcal{B},V^{\mathcal{B}})$.
\end{rema}

\section{Basic cohomology of $\SV$ with trivial coefficients}
In this section, we will compute the basic cohomology groups of $\SV$ with coefficients in its trivial module $\mathbb{C}_a$.
Since ${\rm \tilde
H}^q(\SV,\mathbb{C}_a)\cong {\rm \tilde
H}^q(\SV,\mathbb{C})$ for any $a\in \C$ (\cite{BKV}), we only need to compute ${\rm \tilde
H}^q(\SV,\mathbb{C})$.
In this case, by (\ref{ddd}), the differential $d_n$ of an $n$-cochain $\g$ is
given as follows:
\begin{eqnarray*}
&&(d_n\g)_{\la_1,\cdots,\la_{n+1}}(a_1,\cdots,a_{n+1})\nonumber\\&&\ \ \ =\mbox{$\sum\limits_{1\le i<j\le n+1}$}(-1)^{i+j}\g_{\la_i+\la_j,\la_1,
\cdots,\hat{\la_i},\cdots,\hat{\la_j},\cdots,\la_{n+1}}([a_{i_{\la_i}}a_j],
a_1,\cdots,\hat{a_i},\cdots,\hat{a_j},\cdots,a_{n+1}).
\end{eqnarray*}

\begin{lemm}\label{l0} ${\rm \tilde
H}^0(\SV,\mathbb{C})={\rm H}^0(\SV,\mathbb{C})=\mathbb{C}$.
\end{lemm}
{\it Proof.}
For any $\gamma\in \tilde
C^0(\SV,\mathbb{C})=\mathbb{C}$, $(d_0\gamma)_\lambda (a)=a_\lambda \gamma =0$ for
$a\in \SV$. This means $\tilde D^0(\SV,\mathbb{C})=\mathbb{C}$ and $\tilde B^0(\SV,\mathbb{C})=0$. Thus ${\rm \tilde
H}^0(\SV,\mathbb{C})=\mathbb{C}$. Moreover, $ {\rm H}^0(\SV,\mathbb{C})=\mathbb{C}$
since $\partial\mathbb{C}=0$.
\QED

Let $\gamma\in \tilde
C^q(\SV,\mathbb{C})$ with $q>0$. By Definition \ref{cochain}, $\gamma$ is determined by its value on $a_1\otimes\cdots \otimes a_q$ with $a_i\in\{L, Y, M, N\}$. Since $\gamma$ is skew-symmetric, we can always assume that the first $k$ variables are $L$, the following $l$ variables are $Y$, the $m$ variables are $M$ and the last $n$ variables are $N$ in $\gamma_{\lambda_1,\cdots,\lambda_{q}}(a_1,\cdots,a_q)$.
Thus we can regard $\gamma_{\lambda_1,\cdots,\lambda_{q}}(a_1,\cdots,a_q)$ as a polynomial in $\lambda_1,\cdots,\lambda_q$, which is skew-symmetric in $\lambda_1,\cdots,\lambda_k$,  in $\lambda_{k+1},\cdots,\lambda_{k+l}$, in $\lambda_{k+l+1},\cdots,\lambda_{k+l+m}$, and in $\lambda_{k+l+m+1},\cdots,\lambda_{k+l+m+n}$, respectively, where $q=k+l+m+n$. Therefore, $\gamma_{\lambda_1,\cdots,\lambda_{q}}(a_1,\cdots,a_q)$ is divisible by
$$\mbox{$\prod\limits_{1\leq i< j\leq k}$}(\lambda_i-\lambda_j)\mbox{$\prod\limits_{1\leq i< j\leq l}$}(\lambda_{k+i}-\lambda_{k+j})\mbox{$\prod\limits_{1\leq i< j\leq m}$}(\lambda_{k+l+i}-\lambda_{k+l+j})\mbox{$\prod\limits_{1\leq i< j\leq n}$}(\lambda_{k+l+m+i}-\lambda_{k+l+m+j}),$$
whose degree is $\begin{pmatrix}k\\2\end{pmatrix}+\begin{pmatrix}l\\2\end{pmatrix}+\begin{pmatrix}m\\2\end{pmatrix}+\begin{pmatrix}n\\2\end{pmatrix}$.

Following \cite{BKV}, we define an operator $\tau:\t C^q(\SV,\CC)\rightarrow
\t C^{q-1}(\SV,\CC)$ by
\begin{eqnarray}
(\tau
\g)_{\la_1,\cdots,\la_{q-1}}(a_1,\cdots,a_{q-1})=(-1)^{q-1}\frac{\partial}{\partial\la}\g_{\la_1,\cdots,\la_{q-1},\la}(a_1,\cdots,a_{q-1},L)|_{\la=0}.
\end{eqnarray}
By direct computations (referring to \cite{YW}), we have
\begin{align}\label{al3.6}
((d\tau+&\tau d)\g)_{\lambda_1,\cdots,\lambda_{q}}(a_1,\cdots,a_q)\nonumber\\
=&(-1)^q\frac{\partial}{\partial\lambda}\sum_{i=1}^q(-1)^{i+q+1}\g_{\lambda_i+\lambda,\lambda_1,\cdots,\hat{\lambda_i},\cdots,\lambda_{q}}([{a_i}_{\lambda_i}L],a_1,\cdots,\hat{a_i},\cdots,a_{q})|_{\lambda=0}\nonumber\\
=&\frac{\partial}{\partial\lambda}\sum_{i=1}^q\g_{\lambda_1,\cdots,\lambda_{i-1},\lambda_i+\lambda,\lambda_{i+1},\cdots,\lambda_{q}}(a_1,\cdots,a_{i-1},[{a_i}_{\lambda_i}L],a_{i+1},\cdots,a_{q})|_{\lambda=0}\nonumber\\
=&\frac{\partial}{\partial\lambda}\sum_{i=1}^k(\lambda_i-\lambda)\g_{\lambda_1,\cdots,\lambda_{i-1},\lambda_i+\lambda,\lambda_{i+1},\cdots,\lambda_{q}}(a_1,\cdots,a_{i-1},a_i,a_{i+1},\cdots,a_{q})|_{\lambda=0}\nonumber\\
&+\frac{\partial}{\partial\lambda}\sum_{i=k+1}^{k+l}(\lambda_i-\frac{1}{2}\lambda)\g_{\lambda_1,\cdots,\lambda_{i-1},\lambda_i+\lambda,\lambda_{i+1},\cdots,\lambda_{q}}(a_1,\cdots,a_{i-1},a_i,a_{i+1},\cdots,a_{q})|_{\lambda=0}\nonumber\\
&+\frac{\partial}{\partial\lambda}\sum_{i=k+l+1}^{k+l+m}\lambda_i\g_{\lambda_1,\cdots,\lambda_{i-1},\lambda_i+\lambda,\lambda_{i+1},\cdots,\lambda_{q}}(a_1,\cdots,a_{i-1},a_i,a_{i+1},\cdots,a_{q})|_{\lambda=0}\nonumber\\
=&({\rm deg\,} \g-k-\frac{l}{2})\g_{\lambda_1,\cdots,\lambda_{q}}(a_1,\cdots,a_q),
\end{align}
where ${\rm deg}\, \g$ is the total degree of $\g$ in $\la_1,\cdots,\la_{q}$.
Therefore, if a $q$-cocyle $\g$ satisfies ${\rm deg}\, \g\neq k+\frac{l}{2}$, it must be a coboundary.  Only those homogeneous cochains whose degree as a polynomial is equal to $k+\frac{l}{2}$ contribute to the cohomology of $\t C^\bullet(\SV,\C)$.

Consider the quadratic inequality
\begin{equation}\label{ineq}
\begin{pmatrix}k\\2\end{pmatrix}+\begin{pmatrix}l\\2\end{pmatrix}+\begin{pmatrix}m\\2\end{pmatrix}+\begin{pmatrix}n\\2\end{pmatrix}\leq k+\frac{l}{2}.
\end{equation}

\begin{lemm}\label{le0}
All non-negative integral solutions satisfying $k+\frac{l}{2}\in\Z_+$ of the inequality (\ref{ineq}) are listed in the following table:
\end{lemm}
%%%%%%%%%%%%%%
\begin{CJK*}{GBK}{song}
\setlength{\LTleft}{50pt} \setlength{\LTright}{50pt} %±ížñÓëÒ³Ãæ×óÓÒ±ßÔµÖ®ŒäµÄŸØÀëŸùÎª£°
\begin{longtable}{|c|c|c|c|}\hline
			$q=k+l+m+n$&$(k,l,m,n)$&$\begin{pmatrix}k\\2\end{pmatrix}+\begin{pmatrix}l\\2\end{pmatrix}+\begin{pmatrix}m\\2\end{pmatrix}+\begin{pmatrix}n\\2\end{pmatrix}$&
${\rm deg}\, \g=k+\frac{l}{2}$\\\hline
 0&(0,0,0,0) &0&0\\\hline
     &(1,0,0,0) &0&1\\
  1& (0,0,1,0) &0&0 \\
    &(0,0,0,1)&0&0\\\hline
   &(0,0,1,1)&0&0\\
   &(1,0,1,0)&0&1\\
   2&(1,0,0,1)&0&1\\
   &(2,0,0,0)&1&2\\
   &(0,2,0,0)&1&1\\\hline
    &(1,0,1,1)&0&1\\
    &(1,0,2,0)&1&1\\
    &(1,0,0,2)&1&1\\
    &(2,0,1,0)&1&2\\
    3&(2,0,0,1)&1&2\\
    &(1,2,0,0)&1&2\\
    &(0,2,1,0)&1&1\\
    &(0,2,0,1)&1&1\\
    &(3,0,0,0)&3&3\\\hline
&(1,0,2,1)&1&1\\
&(1,0,1,2)&1&1\\
&(2,0,1,1)&1&2\\
&(2,0,2,0)&2&2\\
&(2,0,0,2)&2&2\\
  4&(3,0,1,0)&3&3\\
  &(3,0,0,1)&3&3\\
  &(0,2,1,1)&1&1\\
  &(1,2,1,0)&1&2\\
  &(1,2,0,1)&1&2\\
  &(2,2,0,0)&2&3\\\hline
   &(3,2,0,0)&4&4\\
   &(3,0,1,1)&3&3\\
   &(2,0,2,1)&2&2\\
   &(2,0,1,2)&2&2\\
   5&(1,2,1,1)&1&2\\
   &(1,2,2,0)&2&2\\
   &(1,2,0,2)&2&2\\
   &(2,2,1,0)&2&3\\
   &(2,2,0,1)&2&3\\\hline
    &(3,2,0,1)&4&4\\
    &(3,2,1,0)&4&4\\
    &(1,2,2,1)&2&2\\
    6&(1,2,1,2)&2&2\\
    &(2,2,1,1)&2&3\\
    &(2,2,2,0)&3&3\\
    &(2,2,0,2)&3&3\\\hline
&(2,2,2,1)&3&3\\
7&(2,2,1,2)&3&3\\
&(3,2,1,1)&4&4\\\hline
	\end{longtable}
\end{CJK*}
	
By the conclusion of the above table, we can obtain the following result immediately.
\begin{lemm}\label{l1}
$\t H^q(\SV,\C)=0,$ if $q\geq8$.
\end{lemm}
\begin{lemm}
$\t H^1(\SV,\C)=0$.
\end{lemm}
{\it Proof.} By the Lemma \ref{le0}, we only need to consider $(k,l,m,n)=(1,0,0,0), (0,0,1,0)$ and $(0,0,0,1)$. Let $\gamma$ be a $1$-cocycle. Inspired by Lemma 3.4 in \cite{WL}, we can deduce that $\gamma_\lambda(L)=\gamma_\lambda(Y)=\gamma_\lambda(M)=0$. Then we can assume that $\gamma_\lambda(N)=a$ for some $a\in\C$. Thus $(d\gamma)_{\lambda_1,\lambda_2}(L,N)=a\lambda_2=0$, implying $\gamma=0$.
\QED

\begin{lemm}\label{l5}
$\t H^2(\SV,\C)=0$.
\end{lemm}
{\it Proof.} We only need to consider $(k,l,m,n)=(0,0,1,1),(1,0,1,0),(1,0,0,1),(2,0,0,0)$ and $(0,2,0,0)$. Let $\gamma$ be a $2$-cocycle. As shown in Lemma 3.5 in \cite{WL}, we can deduce that $\gamma_{\lambda_1,\lambda_2}(L, L)=\gamma_{\lambda_1,\lambda_2}(L, M)=\gamma_{\lambda_1,\lambda_2}(Y,Y)=0$.

Now we can assume that $\gamma_{\lambda_1,\lambda_2}(M, N)=a$ and $\gamma_{\lambda_1,\lambda_2}(L, N)=b\lambda_1+c\lambda_2$, for some $a,b,c\in\C$. By $(d\gamma)_{\lambda_1,\lambda_2, \lambda_3}(L,M,N)=a(\lambda_2+ \lambda_3)=0$ and $(d\gamma)_{\lambda_1,\lambda_2, \lambda_3}(L,L,N)=-b(\lambda_1-\lambda_2)(\lambda_1+\lambda_2+\lambda_3)=0$, we can obtain that $a=b=0$.

Let $\phi$ be a $1$-cochain defined by $\phi_\lambda(N)=c$. It is not difficult to check that $(d\phi)_{\lambda_1,\lambda_2}(L, N)=c\lambda_2$, which implies that $\gamma$ is a coboundary. Replacing $\gamma$ by $\gamma-d\phi$, we can assume that $\gamma_{\lambda_1,\lambda_2}(L, N)=0$.
\QED

\begin{lemm}\label{l6}
$\t H^3(\SV,\C)=\C \Phi^1\oplus\C\Phi^2\oplus\C\Phi^3$, where $\Phi^1,\Phi^2,\Phi^3$ are as follows
\begin{align*}
&{\Phi^1}_{\lambda_1,\lambda_2,\lambda_3}(L, N, N)=\lambda_2-\lambda_3,\\
&{\Phi^2}_{\lambda_1,\lambda_2,\lambda_3}(L, L, N)=(\lambda_1-\lambda_2)\lambda_3,\\
&{\Phi^3}_{\lambda_1,\lambda_2,\lambda_3}(L, L, L)=(\lambda_1-\lambda_2)(\lambda_1-\lambda_3)(\lambda_2-\lambda_3).
\end{align*}
In particular,
${\rm dim}\,\t H^3(\SV,\C)=3$.
\end{lemm}
{\it Proof.} We just need to consider $(k,l,m,n)=(1,0,1,1),(1,0,2,0),(1,0,0,2),(2,0,1,0),(2,0,0,1)$,\\
$(1,2,0,0),(0,2,1,0),(0,2,0,1)$ and $(3,0,0,0)$. Let $\gamma$ be an arbitrary $3$-cocycle. With the similar discussion in Lemma 3.6 in \cite{WL}, we can assume that $\g_{\lambda_1,\lambda_2,\lambda_3}(L, L, M)=\g_{\lambda_1,\lambda_2,\lambda_3}(L, Y, Y)=\g_{\lambda_1,\lambda_2,\lambda_3}(L, M, M)=\g_{\lambda_1,\lambda_2,\lambda_3}(Y,Y,M)=0.$ Then we can assume
\begin{align}\label{LMNbc}
&\gamma_{\lambda_1,\lambda_2,\lambda_3}(L,M,N)=a\lambda_1+b\lambda_2+c\lambda_3,\nonumber\\
&\gamma_{\lambda_1,\lambda_2,\lambda_3}(L,N,N)=e(\lambda_2-\lambda_3),\nonumber\\
&\gamma_{\lambda_1,\lambda_2,\lambda_3}(L,L,N)=(\lambda_1-\lambda_2)(f\lambda_1+f\lambda_2+g\lambda_3),\\
&\gamma_{\lambda_1,\lambda_2,\lambda_3}(Y,Y,N)=h(\lambda_1-\lambda_2),\nonumber\\
&\gamma_{\lambda_1,\lambda_2,\lambda_3}(L, L, L)=k(\lambda_1-\lambda_2)(\lambda_1-\lambda_3)(\lambda_2-\lambda_3),\nonumber
\end{align}
where $a,b,c,e,f,g,h,k\in\C$. In the sequel, we try to determine these parameters.

Let $\varphi$ be a $2$-cochain defined by $\varphi_{\lambda_1,\lambda_2}(M,N)=h$ and $\varphi_{\lambda_1,\lambda_2}(L,N)=f\lambda_1$.
Then
\begin{eqnarray}\label{LMNvarphi}
&&(d\varphi)_{\lambda_1,\lambda_2,\lambda_3}(L,M,N)=h(\lambda_2+\lambda_3),\nonumber\\
&&(d\varphi)_{\lambda_1,\lambda_2,\lambda_3}(Y,Y,N)=-h(\lambda_1-\lambda_2),\\
&&(d\varphi)_{\lambda_1,\lambda_2,\lambda_3}(L,L,N)=-f(\lambda_1-\lambda_2)(\lambda_1+\lambda_2+\lambda_3).\nonumber
\end{eqnarray}
Replacing $\gamma$ by $\gamma+d\varphi$, we can assume $h=f=0$.
Then by direct computations, we have
\begin{eqnarray*}
(d\gamma)_{\lambda_1,\lambda_2,\lambda_3,\lambda_4}(L,Y,Y,N)
=(\lambda_2-\lambda_3)(a\lambda_1+b(\lambda_2+\lambda_3)+c\lambda_4)=0,
\end{eqnarray*}
which implies $a=b=c=0$.

Since $\SV=(\C[\partial]L\oplus\C[\partial]N)\ltimes (\C[\partial]Y\oplus\C[\partial]M)$ and $\C[\partial]L\oplus\C[\partial]N$ is isomorphic to the Heisenberg-Virasoro conformal algebra $\mathcal{HV}$, whose cohomology was studied in \cite{YW}, we can deduce that the cochain defined by
\begin{align*}
&{\Phi^1}_{\lambda_1,\lambda_2,\lambda_3}(L, N, N)=\lambda_2-\lambda_3,\\
&{\Phi^2}_{\lambda_1,\lambda_2,\lambda_3}(L, L, N)=(\lambda_1-\lambda_2)\lambda_3,\\
&{\Phi^3}_{\lambda_1,\lambda_2,\lambda_3}(L, L, L)=(\lambda_1-\lambda_2)(\lambda_1-\lambda_3)(\lambda_2-\lambda_3),
\end{align*}
are cocycles, but not coboundaries. Therefore, we have $\t H^3(\SV,\C)=\C \Phi^1\oplus\C\Phi^2\oplus\C\Phi^3$.
\QED

\begin{rema}
(1) Actually $b$ and $c$ in $\gamma_{\lambda_1,\lambda_2,\lambda_3}(L,M,N)$ in (\ref{LMNbc}) have changed after adding $d\varphi$ to $\g$ in (\ref{LMNvarphi}), but we still use the original notation. For convenience, the following similar issues will not be repeated.\\
(2)The above skew-symmetric function $\Phi^1: \SV\otimes\SV\otimes\SV\to\C[\lambda_1,\lambda_2,\lambda_3]$ has values $\lambda_2-\lambda_3$ on $L\otimes N\otimes N$ and $0$ on others.
 \end{rema}
%Let $\gamma$ be a $3$-cochain defined by  %$\gamma_{\lambda_1,\lambda_2,\lambda_3}(L, M, M)=\lambda_2-\lambda_3$. It %was shown in \cite{BKV} that $\gamma$ is a $3$-cocycle and not a %$3$-coboundary.

\begin{lemm}\label{l8}
$\t H^4(\SV,\C)=\C \Psi^1\oplus\C\Psi^2$, where $\Psi^1,\Psi^2$ are as follows
\begin{align*}
&{\Psi^1}_{\lambda_1,\lambda_2,\lambda_3,\lambda_4}(L, L, N, N)=(\lambda_1-\lambda_2)(\lambda_3-\lambda_4),\\
&{\Psi^2}_{\lambda_1,\lambda_2,\lambda_3,\lambda_4}(L, L, L, N)=(\lambda_1-\lambda_2)(\lambda_1-\lambda_3)(\lambda_2-\lambda_3).
\end{align*}
In particular,
${\rm dim}\,\t H^4(\SV,\C)=2$.
\end{lemm}
{\it Proof.} For $q=4$, we only need to consider $(k,l,m,n)=(1,0,2,1),(1,0,1,2),(2,0,1,1),(2,0,2,0),$\\$(2,0,0,2),(3,0,1,0),(3,0,0,1),(0,2,1,1),(1,2,1,0),(1,2,0,1)$ and $(2,2,0,0)$.
Let $\gamma$ be an arbitrary $4$-cocycle. As mentioned in Lemma 3.8 in \cite{WL}, we can suppose that $\gamma_{\lambda_1,\lambda_2,\lambda_3,\lambda_4}(L, L, M ,M)=\gamma_{\lambda_1,\lambda_2,\lambda_3,\lambda_4}(L, L, L,M)=\gamma_{\lambda_1,\lambda_2,\lambda_3,\lambda_4}(L,Y,Y,M)=\gamma_{\lambda_1,\lambda_2,\lambda_3,\lambda_4}(L,L,Y,Y)=0$.
Then we can assume $\gamma$ is defined by
\begin{align*}
&\gamma_{\lambda_1,\lambda_2,\lambda_3,\lambda_4}(L,M,M,N)=a(\lambda_2-\lambda_3),\\
&\gamma_{\lambda_1,\lambda_2,\lambda_3,\lambda_4}(L,M,N,N)=b(\lambda_3-\lambda_4),\\
&\gamma_{\lambda_1,\lambda_2,\lambda_3,\lambda_4}(L,L,M,N)=(\lambda_1-\lambda_2)(c\lambda_1+c\lambda_2+e\lambda_3+f\lambda_4),\\
&\gamma_{\lambda_1,\lambda_2,\lambda_3,\lambda_4}(L,L,N,N)=g(\lambda_1-\lambda_2)(\lambda_3-\lambda_4),\\
&\gamma_{\lambda_1,\lambda_2,\lambda_3,\lambda_4}(L,L,L,N)=h(\lambda_1-\lambda_2)(\lambda_2-\lambda_3)(\lambda_1-\lambda_3),\\
&\gamma_{\lambda_1,\lambda_2,\lambda_3,\lambda_4}(Y,Y,M,N)=k(\lambda_1-\lambda_2),\\
&\gamma_{\lambda_1,\lambda_2,\lambda_3,\lambda_4}(L,Y,Y,N)=(\lambda_2-\lambda_3)(r\lambda_1+s\lambda_2+s\lambda_3+t\lambda_4),
\end{align*}
where $a,b,c,e,f,g,h,k,r,s,t\in\C$. \\

Let $\phi$ be a $3$-cochain defined by $\phi_{\lambda_1,\lambda_2,\lambda_3}(L,M,N)=\frac{b}{2}\lambda_3$.
Then by direct computations, we have
\begin{align*}
&(d\phi)_{\lambda_1,\lambda_2,\lambda_3,\lambda_4}(L,M,N,N)=b(\lambda_3-\lambda_4).
\end{align*}
Replacing $\gamma$ by $\gamma-d\phi$, we can assume $b=0$. In this case, $c,e,f,r,s,t$ change. It does not matter.  Consequently, $\g_{\lambda_1,\lambda_2,\lambda_3,\lambda_4}(L,M,N,N)=0$.

Let $\psi$ be a $3$-cochain defined by $\psi_{\lambda_1,\lambda_2,\lambda_3}(Y,Y,N)=s(\lambda_1-\lambda_2)$.
Then by direct computations, we have
\begin{align*}
&(d\psi)_{\lambda_1,\lambda_2,\lambda_3,\lambda_4}(L,Y,Y,N)
=s(\lambda_2-\lambda_3)(\lambda_2+\lambda_3+\lambda_4).
\end{align*}
Replacing $\gamma$ by $\gamma-d\psi$, we can assume $s=0$.  Consequently,
$\gamma_{\lambda_1,\lambda_2,\lambda_3,\lambda_4}(L,Y,Y,N)
=(\lambda_2-\lambda_3)(r\lambda_1+t\lambda_4)$.

Since \begin{align*}
&(d\gamma)_{\lambda_1,\lambda_2,\lambda_3,\lambda_4,\lambda_5}(L,Y,Y,M,N)
=(\lambda_2-\lambda_3)((a+k)\lambda_2+(a+k)\lambda_3+(k-a)\lambda_4+k\lambda_5)=0,
\end{align*}
we have $a=k=0$.

Since
\begin{align*}
&(d\gamma)_{\lambda_1,\lambda_2,\lambda_3,\lambda_4,\lambda_5}(L,Y,Y,N,N)
=-2t(\lambda_2-\lambda_3)(\lambda_4-\lambda_5)=0,
\end{align*}
we have $t=0$.

Since \begin{align*}
&(d\gamma)_{\lambda_1,\lambda_2,\lambda_3,\lambda_4,\lambda_5}(L,L,M,N,N)
=2(e-f)(\lambda_1-\lambda_2)(\lambda_4-\lambda_5)=0,
\end{align*}
we have $e=f$.

Since \begin{align*}
&(d\gamma)_{\lambda_1,\lambda_2,\lambda_3,\lambda_4,\lambda_5}(L,L,Y,Y,N)
=-(\lambda_1-\lambda_2)(\lambda_3-\lambda_4)((c+r)(\lambda_1+\lambda_2)+
(e+r)(\lambda_3+\lambda_4)+(f+r)\lambda_5)=0,
\end{align*}
we have $c=e=f=-r$.

Now we can deduce that
 \begin{align*}
 &\g_{\lambda_1,\lambda_2,\lambda_3,\lambda_4}(L,L,M,N)=c(\lambda_1-\lambda_2)(\lambda_1+\lambda_2+\lambda_3+\lambda_4),\\
&\g_{\lambda_1,\lambda_2,\lambda_3,\lambda_4}(L,Y,Y,N)=-c\lambda_1(\lambda_2-\lambda_3).
\end{align*}
Let $\varphi$ be a $3$-cochain defined by $\varphi_{\lambda_1,\lambda_2,\lambda_3}(L,L,M)=\frac{c}{2}(\lambda_1-\lambda_2)(\lambda_1+\lambda_2+\lambda_3)$ and $\varphi_{\lambda_1,\lambda_2,\lambda_3}(L,Y,Y)=-\frac{c}{2}\lambda_1(\lambda_2-\lambda_3)$.
%$\varphi_{\lambda_1,\lambda_2,\lambda_3}(L,M,N)=\varphi_{\lambda_1,\lambda_2,\lambda_3}(Y,Y,N)=0$
Then by direct computations, we have
\begin{align*}
&(d\varphi)_{\lambda_1,\lambda_2,\lambda_3,\lambda_4}(L,L,M,N)=c(\lambda_1-\lambda_2)(\lambda_1+\lambda_2+\lambda_3+\lambda_4),\\
&(d\varphi)_{\lambda_1,\lambda_2,\lambda_3,\lambda_4}(L,Y,Y,N)=-c\lambda_1(\lambda_2-\lambda_3),\\
&(d\varphi)_{\lambda_1,\lambda_2,\lambda_3,\lambda_4}(L,L,L,M)=(d\varphi)_{\lambda_1,\lambda_2,\lambda_3,\lambda_4}(L,L,Y,Y)=0.
%&(d\varphi)_{\lambda_1,\lambda_2,\lambda_3,\lambda_4}(L,M,N,N)=0,\\
\end{align*}
Replacing $\gamma$ by $\gamma-d\varphi$, we have $c=0$. Consequently,
$\g_{\lambda_1,\lambda_2,\lambda_3,\lambda_4}(L,L,M,N)=\g_{\lambda_1,\lambda_2,\lambda_3,\lambda_4}(L,Y,Y,N)=0$.

With the similar discussion in the cohomology theory of the Heisenberg-Virasoro conformal algebra $\mathcal{HV}$ introduced in \cite{YW}, which is isomorphic to $\C[\partial]L\oplus\C[\partial]N$, we can deduce that the cochain defined by
\begin{align*}
&{\Psi^1}_{\lambda_1,\lambda_2,\lambda_3,\lambda_4}(L, L, N, N)=(\lambda_1-\lambda_2)(\lambda_3-\lambda_4),\\
&{\Psi^2}_{\lambda_1,\lambda_2,\lambda_3,\lambda_4}(L, L, L, N)=(\lambda_1-\lambda_2)(\lambda_1-\lambda_3)(\lambda_2-\lambda_3),
\end{align*}
are cocycles, but not coboundaries. Thus, we have $\t H^4(\SV,\C)=\C \Psi^1\oplus\C\Psi^2$.
\QED

%%%%%%%%%%%%%%

\begin{lemm}\label{l9}
$\t H^5(\SV,\C)=0$.
\end{lemm}
{\it Proof.} For $q=5$, we just consider $(k,l,m,n)=(3,2,0,0),(3,0,1,1),(2,0,2,1),(2,0,1,2),(1,2,1,1)$,\\$(1,2,2,0),(1,2,0,2),(2,2,1,0)$ and $(2,2,0,1)$. Let $\gamma$ be an arbitrary $5$-cocycle. As revealed in Lemma 3.9 in \cite{WL}, we can assume that
\begin{align*}
&\gamma_{\lambda_1,\lambda_2,\lambda_3,\lambda_4,\lambda_5}(L,L,L,Y,Y)=0,\\
&\gamma_{\lambda_1,\lambda_2,\lambda_3,\lambda_4, \lambda_5}(L,Y,Y,M,M)=i(\lambda_2-\lambda_3)(\lambda_4-\lambda_5),\\
&\gamma_{\lambda_1,\lambda_2,\lambda_3,\lambda_4,\lambda_5}(L,L,Y,Y,M)=j(\lambda_1-\lambda_2)(\lambda_3-\lambda_4)\lambda_5,
\end{align*}
for some $i,j\in\C.$
Then we can assume that $\gamma$ is defined by
\begin{align*}
&\gamma_{\lambda_1,\lambda_2,\lambda_3,\lambda_4,\lambda_5}(L,L,L,M,N)=a(\lambda_1-\lambda_2)(\lambda_1-\lambda_3)(\lambda_2-\lambda_3),\\
&\gamma_{\lambda_1,\lambda_2,\lambda_3,\lambda_4,\lambda_5}(L,L,M,M,N)=b(\lambda_1-\lambda_2)(\lambda_3-\lambda_4),\\
&\gamma_{\lambda_1,\lambda_2,\lambda_3,\lambda_4, \lambda_5}(L,L,M,N,N)=c(\lambda_1-\lambda_2)(\lambda_4-\lambda_5),\\
&\gamma_{\lambda_1,\lambda_2,\lambda_3,\lambda_4, \lambda_5}(L,Y,Y,M,N)=(\lambda_2-\lambda_3)(e\lambda_1+f\lambda_2+f\lambda_3+g\lambda_4+h\lambda_5),\\
&\gamma_{\lambda_1,\lambda_2,\lambda_3,\lambda_4, \lambda_5}(L,Y,Y,N,N)=k(\lambda_2-\lambda_3)(\lambda_4-\lambda_5),\\
&\gamma_{\lambda_1,\lambda_2,\lambda_3,\lambda_4, \lambda_5}(L,L,Y,Y,N)=(\lambda_1-\lambda_2)(\lambda_3-\lambda_4)(x\lambda_1+x\lambda_2+y\lambda_3+y\lambda_4+z\lambda_5),
\end{align*}
where $a,b,c,e,f,g,h,k,x,y,z\in\C$.

Let $\phi$ be a $4$-cochain defined by
$\phi_{\lambda_1,\lambda_2,\lambda_3,\lambda_4}(L,L,M,N)=\frac{z-y}{2}(\lambda_1-\lambda_2)(\lambda_3-\lambda_4),$ $\phi_{\lambda_1,\lambda_2,\lambda_3,\lambda_4}(L,L,Y,Y)=-(\lambda_1-\lambda_2)(\lambda_3-\lambda_4)(\frac{x}{2}\lambda_1+\frac{x}{2}\lambda_2+\frac{y+z}{4}\lambda_3+\frac{y+z}{4}\lambda_4)$ and $\phi_{\lambda_1,\lambda_2,\lambda_3,\lambda_4}(L,Y,Y,N)=0$.
Then $d\phi$ is given by
\begin{align*}
&(d\phi)_{\lambda_1,\lambda_2,\lambda_3,\lambda_4,\lambda_5}(L,L,Y,Y,N)=(\lambda_1-\lambda_2)(\lambda_3-\lambda_4)(x\lambda_1+x\lambda_2+y\lambda_3+y\lambda_4+z\lambda_5),\\
&(d\phi)_{\lambda_1,\lambda_2,\lambda_3,\lambda_4,\lambda_5}(L,L,L,Y,Y)=0.
\end{align*}
Replacing $\gamma$ by $\gamma-d\phi$, we have
$\gamma_{\lambda_1,\lambda_2,\lambda_3,\lambda_4,\lambda_5}(L,L,Y,Y,N)=0$.

Define
\begin{eqnarray*}
&&\psi_{\lambda_1,\lambda_2,\lambda_3,\lambda_4}(Y,Y,M,N)=f(\lambda_1-\lambda_2),\\
&&\psi_{\lambda_1,\lambda_2,\lambda_3,\lambda_4}(L,Y,Y,M)=(\lambda_2-\lambda_3)(-\frac{e}{4}\lambda_1+\frac{f-g}{4}\lambda_4),\\
&&\psi_{\lambda_1,\lambda_2,\lambda_3,\lambda_4}(L,M,M,N)=0.
\end{eqnarray*}
Then we can deduce that
 \begin{eqnarray*}
(d\psi)_{\lambda_1,\lambda_2,\lambda_3,\lambda_4, \lambda_5}(L,Y,Y,M,N)=(\lambda_2-\lambda_3)(e\lambda_1+f\lambda_2+f\lambda_3+g\lambda_4+\frac{f+g}{2}\lambda_5.
\end{eqnarray*}
Replacing $\gamma$ by $\gamma-d\psi$, we can assume that $$\gamma_{\lambda_1,\lambda_2,\lambda_3,\lambda_4,\lambda_5}(L,Y,Y,M,N)=h(\lambda_2-\lambda_3)\lambda_5$$ for some $h\in\C$.

It is not difficult to check that
\begin{align*}
&(d\gamma)_{\lambda_1,\lambda_2,\lambda_3,\lambda_4,\lambda_5,\lambda_6}(L,L,Y,Y,M,N)=(\lambda_1-\lambda_2)(\lambda_3-\lambda_4)(-b\lambda_3-b\lambda_4+(4j+b)\lambda_5+2j\lambda_6)=0,\\
&(d\gamma)_{\lambda_1,\lambda_2,\lambda_3,\lambda_4,\lambda_5,\lambda_6}(L,Y,Y,M,M,N)=6i(\lambda_2-\lambda_3)(\lambda_4-\lambda_5)=0,\\
&(d\gamma)_{\lambda_1,\lambda_2,\lambda_3,\lambda_4,\lambda_5,\lambda_6}(L,Y,Y,M,N,N)=4h(\lambda_2-\lambda_3)(\lambda_5-\lambda_6)=0,
\end{align*}
which implies $b=j=i=h=0$. Thus, $\gamma_{\lambda_1,\lambda_2,\lambda_3,\lambda_4,\lambda_5}(L,L,M,M,N)=\gamma_{\lambda_1,\lambda_2,\lambda_3,\lambda_4,\lambda_5}(L,L,Y,Y,M)=\gamma_{\lambda_1,\lambda_2,\lambda_3,\lambda_4,\lambda_5}(L,Y,Y,M,M)=\gamma_{\lambda_1,\lambda_2,\lambda_3,\lambda_4,\lambda_5}(L,Y,Y,M,N)=0$.
%%%%%%%%%%%%%%%

Similarly, define
\begin{equation*}
\varphi_{\lambda_1,\lambda_2,\lambda_3,\lambda_4}(L,M,N,N)=k(\lambda_3-\lambda_4),
\end{equation*}
then
\begin{equation*}
(d\varphi)_{\lambda_1,\lambda_2,\lambda_3,\lambda_4, \lambda_5}(L,Y,Y,N,N)=k(\lambda_2-\lambda_3)(\lambda_4-\lambda_5).
\end{equation*}Replacing $\gamma$ by $\gamma-d\varphi$, we have $\gamma_{\lambda_1,\lambda_2,\lambda_3,\lambda_4,\lambda_5}(L,Y,Y,N,N)=0$.

Then by direct computations, we have
\begin{align*}
&{(d\g)}_{\lambda_1,\lambda_2,\lambda_3,\lambda_4,\lambda_5,\lambda_6}(L,L,L,Y,Y,N)=a(\lambda_1-\lambda_2)(\lambda_1-\lambda_3)(\lambda_2-\lambda_3)(\lambda_4-\lambda_5)=0,\\
&{(d\g)}_{\lambda_1,\lambda_2,\lambda_3,\lambda_4,\lambda_5,\lambda_6}(L,L,Y,Y,N,N)=-c(\lambda_1-\lambda_2)(\lambda_3-\lambda_4)(\lambda_5-\lambda_6)=0,
\end{align*}
which implying $a=c=0$. Thus, we can deduce that $\t H^5(\SV,\C)=0$.
\QED

\begin{lemm}\label{l10}
$\t H^6(\SV,\C)=0$.
\end{lemm}
{\it Proof.} For $q=6$, we only need to consider $(k,l,m,n)=(3,2,0,1),(3,2,1,0),(1,2,2,1),(1,2,1,2)$,\\$(2,2,1,1),(2,2,2,0)$ and $(2,2,0,2)$. Let $\gamma$ be an arbitrary $6$-cocycle, we can assume that $\gamma$ satisfies
\begin{align}\label{a3.4}
&\gamma_{\lambda_1,\lambda_2,\lambda_3,\lambda_4,\lambda_5,\lambda_6}(L,L,L,Y,Y,N)=a(\lambda_1-\lambda_2)(\lambda_1-\lambda_3)
(\lambda_2-\lambda_3)(\lambda_4-\lambda_5),\nonumber\\
&\gamma_{\lambda_1,\lambda_2,\lambda_3,\lambda_4,\lambda_5,\lambda_6}(L,L,L,Y,Y,M)=b(\lambda_1-\lambda_2)(\lambda_1-\lambda_3)
(\lambda_2-\lambda_3)(\lambda_4-\lambda_5),\nonumber\\
&\gamma_{\lambda_1,\lambda_2,\lambda_3,\lambda_4,\lambda_5,\lambda_6}(L,Y,Y,M,M,N)=c(\lambda_2-\lambda_3)(\lambda_4-\lambda_5),\nonumber\\
&\gamma_{\lambda_1,\lambda_2,\lambda_3,\lambda_4,\lambda_5,\lambda_6}(L,Y,Y,M,N,N)=e(\lambda_2-\lambda_3)(\lambda_5-\lambda_6),\\
&\gamma_{\lambda_1,\lambda_2,\lambda_3,\lambda_4,\lambda_5,\lambda_6}(L,L,Y,Y,M,N)=(\lambda_1-\lambda_2)(\lambda_3-\lambda_4)(f\lambda_1+f\lambda_2+g\lambda_3+g\lambda_4+h\lambda_5+i\lambda_6),\nonumber\\
&\gamma_{\lambda_1,\lambda_2,\lambda_3,\lambda_4,\lambda_5,\lambda_6}(L,L,Y,Y,M,M)=j(\lambda_1-\lambda_2)(\lambda_3-\lambda_4)(\lambda_5-\lambda_6),\nonumber\\
&\gamma_{\lambda_1,\lambda_2,\lambda_3,\lambda_4,\lambda_5,\lambda_6}(L,L,Y,Y,N,N)=k(\lambda_1-\lambda_2)(\lambda_3-\lambda_4)(\lambda_5-\lambda_6),\nonumber
\end{align}
for some $a,b,c,e,f,g,h,i,j,k\in\C$.

Firstly, let $\varphi$ be a $5$-cochain defined by
\begin{align*}
&\varphi_{\lambda_1,\lambda_2,\lambda_3,\lambda_4,\lambda_5}(L,Y,Y,M,N)=\frac{e}{4}(\lambda_2-\lambda_3)\lambda_5,\\
&\varphi_{\lambda_1,\lambda_2,\lambda_3,\lambda_4,\lambda_5}(L,L,Y,Y,M)=\frac{1}{4}(\lambda_1-\lambda_2)(\lambda_3-\lambda_4)(f\lambda_1+f\lambda_2+g\lambda_3+g\lambda_4+h\lambda_5),\\
&\varphi_{\lambda_1,\lambda_2,\lambda_3,\lambda_4,\lambda_5}(L,Y,Y,M,M)=\frac{c}{6}(\lambda_2-\lambda_3)(\lambda_4-\lambda_5),\\
&\varphi_{\lambda_1,\lambda_2,\lambda_3,\lambda_4,\lambda_5}(L,L,M,N,N)=-k(\lambda_1-\lambda_2)(\lambda_4-\lambda_5),\\
&\varphi_{\lambda_1,\lambda_2,\lambda_3,\lambda_4,\lambda_5}(L,L,L,M,N)=a(\lambda_1-\lambda_2)(\lambda_1-\lambda_3)(\lambda_2-\lambda_3).
\end{align*}
It is not difficult to check that
\begin{align*}
&(d\varphi)_{\lambda_1,\lambda_2,\lambda_3,\lambda_4,\lambda_5,\lambda_6}(L,L,Y,Y,M,N)=(\lambda_1-\lambda_2)(\lambda_3-\lambda_4)(f\lambda_1+f\lambda_2+g\lambda_3+g\lambda_4+h\lambda_5+\frac{g+h}{2}\lambda_6),\\
&(d\varphi)_{\lambda_1,\lambda_2,\lambda_3,\lambda_4,\lambda_5,\lambda_6}(L,Y,Y,M,M,N)=c(\lambda_2-\lambda_3)(\lambda_4-\lambda_5),\\
&(d\varphi)_{\lambda_1,\lambda_2,\lambda_3,\lambda_4,\lambda_5,\lambda_6}(L,Y,Y,M,N,N)=e(\lambda_2-\lambda_3)(\lambda_5-\lambda_6),\\
&(d\varphi)_{\lambda_1,\lambda_2,\lambda_3,\lambda_4,\lambda_5,\lambda_6}(L,L,Y,Y,N,N)=k(\lambda_1-\lambda_2)(\lambda_3-\lambda_4)(\lambda_5-\lambda_6),\\
&(d\varphi)_{\lambda_1,\lambda_2,\lambda_3,\lambda_4,\lambda_5,\lambda_6}(L,L,L,Y,Y,N)=a(\lambda_1-\lambda_2)(\lambda_1-\lambda_3)(\lambda_2-\lambda_3) (\lambda_4-\lambda_5).
\end{align*}
Replacing $\gamma$ by $\gamma-d\varphi$, we can assume $a=c=e=f=g=h=k=0$.  Consequently, we can deduce that all the nonzero equations in (\ref{a3.4}) are as follows
\begin{align*}
&\gamma_{\lambda_1,\lambda_2,\lambda_3,\lambda_4,\lambda_5,\lambda_6}(L,L,Y,Y,M,N)=i(\lambda_1-\lambda_2)(\lambda_3-\lambda_4)\lambda_6,\\
&\gamma_{\lambda_1,\lambda_2,\lambda_3,\lambda_4,\lambda_5,\lambda_6}(L,L,L,Y,Y,M)=b(\lambda_1-\lambda_2)(\lambda_1-\lambda_3)
(\lambda_2-\lambda_3)(\lambda_4-\lambda_5),\\
&\gamma_{\lambda_1,\lambda_2,\lambda_3,\lambda_4,\lambda_5,\lambda_6}(L,L,Y,Y,M,M)=j(\lambda_1-\lambda_2)(\lambda_3-\lambda_4)(\lambda_5-\lambda_6).
\end{align*}

Secondly, one can check that
\begin{align*}
&(d\gamma)_{\lambda_1,\lambda_2,\lambda_3,\lambda_4,\lambda_5,\lambda_6,\lambda_7}(L,L,Y,Y,M,N,N)=-4i(\lambda_1-\lambda_2)(\lambda_3-\lambda_4)(\lambda_6-\lambda_7)=0,\\
&(d\gamma)_{\lambda_1,\lambda_2,\lambda_3,\lambda_4,\lambda_5,\lambda_6,\lambda_7}(L,L,Y,Y,M,M,N)=-6j(\lambda_1-\lambda_2)(\lambda_3-\lambda_4)(\lambda_5-\lambda_6)=0,\\
&(d\gamma)_{\lambda_1,\lambda_2,\lambda_3,\lambda_4,\lambda_5,\lambda_6,\lambda_7}(L,L,L,Y,Y,M,N)=-4b(\lambda_1-\lambda_2)(\lambda_1-\lambda_3)(\lambda_2-\lambda_3)(\lambda_4-\lambda_5)=0,
\end{align*}
which imply $i=j=b=0$.
Thus, we can deduce that $\t H^6(\SV,\C)=0$.
\QED

\begin{lemm}\label{l11}
$\t H^7(\SV,\C)=0$.
\end{lemm}
{\it Proof.}
For $q=7$, we only need to consider $(k,l,m,n)=(2,2,2,1),(2,2,1,2)$ and $(3,2,1,1)$. Let $\gamma$ be an arbitrary $7$-cocycle, we can assume that $\gamma$ satisfies
\begin{align*}
&\gamma_{\lambda_1,\lambda_2,\lambda_3,\lambda_4,\lambda_5,\lambda_6,\lambda_7}(L,L,Y,Y,M,M,N)=a(\lambda_1-\lambda_2)(\lambda_3-\lambda_4)
(\lambda_5-\lambda_6),\\
&\gamma_{\lambda_1,\lambda_2,\lambda_3,\lambda_4,\lambda_5,\lambda_6,\lambda_7}(L,L,Y,Y,M,N,N)=b(\lambda_1-\lambda_2)(\lambda_3-\lambda_4)
(\lambda_6-\lambda_7),\\
&\gamma_{\lambda_1,\lambda_2,\lambda_3,\lambda_4,\lambda_5,\lambda_6,\lambda_7}(L,L,L,Y,Y,M,N)=c(\lambda_1-\lambda_2)(\lambda_1-\lambda_3)
(\lambda_2-\lambda_3)(\lambda_4-\lambda_5),
\end{align*}
for some $a,b,c\in\C$.

By the proof of Lemma \ref{l10}, $\gamma=d \varphi$ where
$\varphi$ is a $6$-cochain defined by
\begin{align*}
&\varphi_{\lambda_1,\lambda_2,\lambda_3,\lambda_4,\lambda_5,\lambda_6}(L,L,Y,Y,M,M)=-\frac{a}{6}(\lambda_1-\lambda_2)(\lambda_3-\lambda_4)
(\lambda_5-\lambda_6),\\
&\varphi_{\lambda_1,\lambda_2,\lambda_3,\lambda_4,\lambda_5,\lambda_6}(L,L,Y,Y,M,N)=-\frac{b}{4}(\lambda_1-\lambda_2)(\lambda_3-\lambda_4)
\lambda_6,\\
&\phi_{\lambda_1,\lambda_2,\lambda_3,\lambda_4,\lambda_5,\lambda_6}(L,L,L,Y,Y,M)=-\frac{c}{4}(\lambda_1-\lambda_2)(\lambda_1-\lambda_3)
(\lambda_2-\lambda_3)(\lambda_4-\lambda_5).
\end{align*} Thus, $\t H^7(\SV,\C)=0$.
\QED

 As mentioned above, we can obtain the main result of this section as follows.
\begin{theo}\label{th1}The dimension of ${\t H}^q (\SV,\C)$ is given by
\begin{eqnarray*}
{\rm dim\,\t H}^q(\SV,\C)=\left\{
\begin{array}{ll}
1 &{\mbox if}\ q=0,\\
3 &{\mbox if}\ q=3,\\
2 &{\mbox if}\ q=4,\\
0 &{\mbox otherwise}.
\end{array}
\right.
\end{eqnarray*}
In particular,
\begin{eqnarray}\label{hd001}
{\t H^q(\SV, \C)}=\left\{
\begin{array}{ll}
\C  &{\mbox if}\ q=0,\\
\C \Phi^1\oplus\C\Phi^2\oplus\C\Phi^3 &{\mbox if}\ q=3,\\
\C \Psi^1\oplus\C\Psi^2 &{\mbox if}\ q=4,\\
0 &{\mbox otherwise},
\end{array}
\right.
\end{eqnarray}
where
\begin{align*}
&{\Phi^1}_{\lambda_1,\lambda_2,\lambda_3}(L, N, N)=\lambda_2-\lambda_3,\\
&{\Phi^2}_{\lambda_1,\lambda_2,\lambda_3}(L, L, N)=(\lambda_1-\lambda_2)\lambda_3,\\
&{\Phi^3}_{\lambda_1,\lambda_2,\lambda_3}(L, L, L)=(\lambda_1-\lambda_2)(\lambda_1-\lambda_3)(\lambda_2-\lambda_3),\\
&{\Psi^1}_{\lambda_1,\lambda_2,\lambda_3,\lambda_4}(L, L, N, N)=(\lambda_1-\lambda_2)(\lambda_3-\lambda_4),\\
&{\Psi^2}_{\lambda_1,\lambda_2,\lambda_3,\lambda_4}(L, L, L, N)=(\lambda_1-\lambda_2)(\lambda_1-\lambda_3)(\lambda_2-\lambda_3).
\end{align*}
\end{theo}
{\it Proof.} It follows directly from the previous Lemmas.
\QED

\begin{rema}
The corresponding annihilation algebra of $\SV$ is
\begin{eqnarray*}
&&\textit{Lie}(\SV)^+= \sum_{m\geq -1}\CC L_m+\sum_{n\geq 0}\CC M_n +\sum_{k\geq 0}\CC N_k +\sum_{p\in\frac{1}{2}+\Z_+}\CC Y_{p},\\
&&[L_m,L_{n}]=(m-n)L_{m+n},\ \ \
[L_m,M_n]=-nM_{m+n},\\
&&[\,Y_p\,,Y_{q}\,]=(p-q)M_{p+q},\ \ \ \ \
[\,L_m,Y_p\,]=(\frac{m}{2}-p)Y_{m+p},\\
&&[L_m,N_k]=-kN_{m+k},\ \ \
[N_k,Y_p]=Y_{k+p},\ \ [N_k,M_n]=2M_{k+n},
\end{eqnarray*}
which is a 'half part' of the extended Schr\"odinger-Virasoro Lie algebra (\cite{SY}). So by \cite{BKV}, the dimension of all the cohomology groups of $\textit{Lie}(\SV)^+$ is given by
\begin{eqnarray*}
{\rm dim\, H}^q(\textit{Lie}(\SV)^+,\C)=\left\{
\begin{array}{ll}
1 &{\mbox if}\ q=0,\\
3 &{\mbox if}\ q=3,\\
2 &{\mbox if}\ q=4,\\
0 &{\mbox otherwise}.
\end{array}
\right.
\end{eqnarray*}
 \end{rema}

\section{Reduced cohomology of $\SV$ with trivial coefficients}
In this section, we compute the reduced cohomology groups of $\SV$ with coefficients in its trivial module $\mathbb{C}_a$.

\subsection{Computation of $H^q (\SV,\C)$}
\begin{theo}\label{th2} The dimension of $H^q (\SV,\C)$ is given by
\begin{eqnarray*}
{\rm dim\, H}^q(\SV,\C)=\left\{
\begin{array}{ll}
1 &{\mbox if}\ q=0,\\
3 &{\mbox if}\ q=2,\\
5 &{\mbox if}\ q=3,\\
2 &{\mbox if}\ q=4,\\
0 &{\mbox otherwise}.
\end{array}
\right.
\end{eqnarray*}
\end{theo}
{\it Proof.}
By Proposition \ref{pro2.12}, the map $\g \mapsto \partial \cdot \g$ gives an isomorphism such that $\t H^q(SV, \C)\cong H^q(\partial\t C^\bullet) $ for all $q\geq1$. Therefore, we can obtain the following result immediately by the discussion of Section 3.

\begin{eqnarray}\label{hd}
{H^q(\partial\t C^\bullet) }=\left\{
\begin{array}{llllll}
\C (\partial\Phi^1)\oplus\C(\partial\Phi^2)\oplus\C(\partial\Phi^3) &{\mbox if}\ q=3,\\
\C (\partial\Psi^1)\oplus\C(\partial\Psi^2) &{\mbox if}\ q=4,\\
0 &{\mbox otherwise}.
\end{array}
\right.
\end{eqnarray}

Similar to the discussions in Section 2, we can obtain the following long exact sequence of cohomology groups:
\begin{align}\label{longSVC}
\cdots\longrightarrow& \ \ H^q(\partial\t C^\bullet) \   \stackrel{i_q}{\longrightarrow} \ \ \t H^q(\SV,\C)\  \stackrel{p_q}{\longrightarrow} \ \ H^q (\SV,\C) \ \stackrel{w_q}\longrightarrow\\
\longrightarrow& H^{q+1}(\partial\t C^\bullet) \stackrel{i_{q+1}}{\longrightarrow} \t H^{q+1}(\SV,\C) \stackrel{p_{q+1}}{\longrightarrow} H^{q+1} (\SV,\C) \longrightarrow\cdots\nonumber
\end{align}
where $i_q, p_q$ are induced by $i,p$ in (\ref{longexact}) respectively and $w_q$ is the $q$-th connecting homomorphism. Take $\partial \g \in H^q(\partial\t C^\bullet) $ with a nonzero element $\g\in\t H^q(\SV,\C)$ of degree $k+\frac{l}{2}$, we can obtain that $i_q(\partial \g)=\partial \g\in \t H^q(\SV,\C)$.
Since $\partial\gamma=(\sum\lambda_i)\gamma$, we have $${\rm deg}\,(\partial \g)={\rm deg}\,(\g)+1=k+\frac{l}{2}+1.$$
Thus, $\partial \g=0 \in\t H^q(\SV,\C)$ (Here we use the property $a=0$). Thus, the image of $i_q$ is zero. Then the long exact sequence (\ref{longSVC}) splits into
the following short exact sequence immediately:
\begin{align}\label{iq}
0\longrightarrow \t H^q(\SV,\C)\stackrel{p_q}{\longrightarrow} H^q (\SV,\C)\stackrel{w_q}\longrightarrow H^{q+1}(\partial\t C^\bullet)\longrightarrow 0,
\end{align}
for all $q\geq1$. Thus, we have
\begin{align}\label{dimc}
{\rm dim}\, H^q (\SV,\C) =& {\rm dim}\, \t H^q(\SV,\C)+ {\rm dim}\, H^{q+1}(\partial\t C^\bullet)\nonumber \\
=&{\rm dim}\,\t H^q(\SV,\C)+{\rm dim}\, \t H^{q+1}(\SV,\C)
\end{align}
for all $q\geq1$.
This completes the proof.
\QED

\begin{rema}It was shown in \cite{SY} that there is a unique nontrivial extension of $\SV$ by a $3$-dimensional center.
This coincides with our result ${\rm dim\, H}^2(\SV,\C)=3$.
\end{rema}

It is not difficult to check by (\ref{iq}) that the basis of $H^q(\SV,\C)$ can be obtained by combining the images of a basis of $\t H^q(\SV,\C)$ with the pre-images of a basis of $\t H^{q+1}(\SV,\C)$. Let $\g$ be a nonzero $(q+1)$-cocycle of degree $k+\frac{l}{2}$ such that $\partial\g \in H^{q+1}(\partial\t C^\bullet)$. By (\ref{al3.6}), we can obtain that
\begin{align}\label{dimcc}
d(\tau(\partial\g))=(d\tau+\tau d)(\partial\g)=(deg(\partial\g)-k-\frac{l}{2})(\partial\g)=((k+\frac{l}{2}+1)-k-\frac{l}{2})(\partial\g)=\partial\g.
\end{align}
Thus, the pre-image $w_q^{-1}(\partial\g)$ of $\partial\g$ under the connecting homomorphism $w_q$ is $\tau(\partial\g)$, i.e., $w_q^{-1}(\partial\g)=\tau(\partial\g).$

Due to the above discussions, we can obtain a specific system to calculate the basis of $H^q(\SV,\C)$ so that we can describe the structure of $H^q(\SV,\C)$ more clearly.

\begin{coro}\label{cor4.3}
\begin{eqnarray}\label{hd2}
{ H^q(\SV, \C)}=\left\{
\begin{array}{ll}
\C &{\mbox if}\ q=0,\\
\C \widehat{\Phi^1}\oplus\C\widehat{\Phi^2}\oplus\C\widehat{\Phi^3} &{\mbox if}\ q=2,\\
\C \Phi^1\oplus\C\Phi^2\oplus\C\Phi^3\oplus\C\widehat{\Psi^1}\oplus\C\widehat{\Psi^2} &{\mbox if}\ q=3,\\
\C \Psi^1\oplus\C\Psi^2 &{\mbox if}\ q=4,\\
0 &{\mbox otherwise},
\end{array}
\right.
\end{eqnarray}
where $\widehat{X}=\tau(\partial X)$ and $X\in \{ \Phi^1,\Phi^2,\Phi^3,\Psi^1,\Psi^2\}$ as shown in Theorem \ref{th1}. More specifically,
\begin{align*}
\widehat{\Phi^1}&_{\lambda_1,\lambda_2}(N,N)=\lambda_1-\lambda_2,\\
\widehat{\Phi^2}&_{\lambda_1,\lambda_2}(L,N)=\lambda_2^2,\\
\widehat{\Phi^3}&_{\lambda_1,\lambda_2}(L,L)=-\lambda_1^3+\lambda_2^3,\\
\widehat{\Psi^1}&_{\lambda_1,\lambda_2,\lambda_3}(L,N,N)=\lambda_2^2-\lambda_3^2,\\
\widehat{\Psi^2}&_{\lambda_1,\lambda_2,\lambda_3}(L,L,N)=-\lambda_1^3-\lambda_1^2\lambda_3+\lambda_2^3+\lambda_2^2\lambda_3.
\end{align*}
\end{coro}
{\it Proof.} It follows directly from the Theorem \ref{th1}, (\ref{dimc}) and (\ref{dimcc}). About the concrete expression of $\widehat{\Phi^1},\widehat{\Phi^2},\widehat{\Phi^3},\widehat{\Psi^1},\widehat{\Psi^2}$, we just take $\widehat{\Phi^1}$ for example, others can be proved similarly.
\begin{align*}
\widehat{\Phi^1}_{\lambda_1,\lambda_2}(N,N)&=(\tau(\partial\Phi^1))_{\lambda_1,\lambda_2}(N,N)\\
&=(-1)^2\frac{\partial}{\partial\lambda}(\partial\Phi^1)_{\lambda_1,\lambda_2,\lambda}(N,N,L)|_{\lambda=0}\\
&=\frac{\partial}{\partial\lambda}(\lambda_1+\lambda_2+\lambda)\Phi^1_{\lambda_1,\lambda_2,\lambda}(N,N,L)|_{\lambda=0}\\
&=\frac{\partial}{\partial\lambda}(\lambda_1+\lambda_2+\lambda)(\lambda_1-\lambda_2)|_{\lambda=0}\\
&=\lambda_1-\lambda_2.
\end{align*}
\QED

\subsection{Computation of $H^q (\SV,\C_a)$ if $a\neq0$}
\begin{theo} \label{th3}
For any $q\in\Z_+$,
${H}^q(\SV,\C_a)=0$ if $a\neq 0$.
\end{theo}
\noindent{\it Proof.~} Similar to the proof of Lemma 3.2 in \cite{YW}, we can define an operator $\tau_2:\t C^q(\SV,\C_a)\rightarrow \t C^{q-1}(\SV,\C_a)$ by
\begin{eqnarray}\label{7-3}
(\tau_2
\g)_{\la_1,\cdots,\la_{q-1}}(a_1,\cdots,a_{q-1})
=(-1)^{q-1}\g_{\la_1,\cdots,\la_{q-1},\la}(a_1,\cdots,a_{q-1},L)|_{\la=0}.
\end{eqnarray}
Then
\begin{align*}
(d\tau_2+\tau_2 d)
\g\equiv -a \g \ (\mbox{mod}\
\partial\t C^q(\SV, \CC_a)),
\end{align*}
which implies ${H}^q(\SV,\C_a)=0$ for all $q\geq0$ if $a\neq 0$.
\QED

\section{Relative cohomology of $\SV$ with trivial coefficients}
In this section, we compute the basic relative cohomology and reduced relative cohomology of $\SV$ modulo $\mathcal{B}$ with coefficients in a trivial module $\mathbb{C}_a$, where $\mathcal{B}$ is a subalgebra of $\SV$. Then we can reformulate our main result using the language of relative cohomology.

Similar to Lemma \ref{l0}, we have
\begin{lemm}
${\rm \tilde
H}^0(\SV;\mathcal{B}, \mathbb{C})={\rm H}^0(\SV;\mathcal{B}, \mathbb{C})=\mathbb{C}$.
\end{lemm}

When we compute the basic relative cohomology, we also need the following lemma.
\begin{lemm}
$\tau (\t C^q(\SV;\mathcal{B},\CC))\subset\t C^q(\SV;\mathcal{B},\CC)$.
\end{lemm}
\noindent{\it Proof.~}
Recall the operator $\tau:\t C^q(\SV,\CC)\rightarrow
\t C^{q-1}(\SV,\CC)$ is given by
\begin{eqnarray*}
(\tau
\g)_{\la_1,\cdots,\la_{q-1}}(a_1,\cdots,a_{q-1})=
(-1)^{q-1}\frac{\partial}{\partial\la}
\g_{\la_1,\cdots,\la_{q-1},\la}(a_1,\cdots,a_{q-1},L)|_{\la=0}.
\end{eqnarray*}
If $\g\in \t C^q(\SV;\mathcal{B},\CC)$, then $\g_{\la_1,\cdots,\la_q}(a_1,\cdots, a_q)=(d\g)_{\la_1,\cdots,\la_{q+1}}(b_1,\cdots, b_{q+1})=0$ for any $a_1,b_1\in\mathcal{B}$. Consequently, $\g_{\la_1,\cdots,\la_{q-1},\la}(a_1,\cdots,a_{q-1},L)|_{\la=0}=0$ and then $(\tau
\g)_{\la_1,\cdots,\la_{q-1}}(a_1,\cdots,a_{q-1})=0$ for any $a_1\in\mathcal{B}$.
So $\tau\g\in \t C^q(\SV;\mathcal{B},\CC)$.
Since $(d\tau+\tau d)\g=({\rm deg\,} \g-k-\frac{l}{2})\g$, then $$d(\tau \g)=({\rm deg\,} \g-k-\frac{l}{2})\g-\tau (d\g)\in \t C^q(\SV;\mathcal{B},\CC).$$
\QED

By the above lemma, the technique of ``comparing degree", developed in the front of Section 3, also works in the computation of basic relative cohomology. Then we can compute the relative cohomology, similarly.

When we compute the reduced relative cohomology ${\rm \t H}^q (\SV;\mathcal{B},\C_a)$ in case $a\neq0$, we need the following lemma.
\begin{lemm}
$\tau_2 (\t C^q(\SV;\mathcal{B},\CC))\subset\t C^q(\SV;\mathcal{B},\CC)$.
\end{lemm}
\noindent{\it Proof.~} It follows from the equation \ref{7-3} immediately.
\QED

Then we have
\begin{theo}
For any $q\in\Z_+$,
${H}^q(\SV;\mathcal{B},\C_a)=0$ if $a\neq 0$.
\end{theo}

Since ${\rm \tilde
H}^q(\SV;\mathcal{B},\mathbb{C}_a)\cong {\rm \tilde
H}^q(\SV;\mathcal{B},\mathbb{C})$ for any $a\in \C$, we only need to consider ${\rm \tilde
H}^q(\SV;\mathcal{B},\mathbb{C})$ and ${\rm
H}^q(\SV;\mathcal{B},\mathbb{C})$ for $q>0$.
In the sequel, we use the notation $\langle X_1, X_2,\cdots, X_n\rangle$ to denote $\C[\partial]X_1+\C[\partial]X_2+\cdots+ \C[\partial]X_n$, where $X_i\in\SV,\ 1\le i\le n$. Then we compute the relative cohomology of $\SV$ module $\mathcal{B}$ with trivial coefficients, where $\mathcal{B}=\langle X_1, X_2,\cdots, X_n\rangle, X_i=L,N,Y,M$.

%\subsection{}
\begin{theo} If $\mathcal{B}=\langle L\rangle, \langle L,N\rangle,
\langle L,M\rangle, \langle L, N,M\rangle$ or $\langle L,Y,M\rangle$, then
${\rm \t H}^q (\SV;\mathcal{B},\C)=0$ and
${\rm H}^q (\SV;\mathcal{B},\C)=0$ for $q>0$.
\end{theo}
\noindent{\it Proof.~} We only consider $\mathcal{B}=\langle L\rangle$, the others are similar. In this case, we need to consider all solutions in Lemma \ref{le0} with $k=0$. Then by Lemma \ref{le0}, $$(l,m,n)=(0,1,0), (0,0,1),(0,1,1),(2,0,0),(2,1,0),(2,0,1),(2,1,1).$$

(i) For $q=1$, assume $\gamma\in {\t C}^1 (\SV;\mathcal{B},\C)$ and denote
$\gamma_\lambda(M)=a, \gamma_\lambda(N)=b$. Since $(d\gamma)_{\lambda_1,\lambda_2}(L,M)=a\lambda_2=0$ and $(d\gamma)_{\lambda_1,\lambda_2}(L,N)=b\lambda_2=0$, we have $a=b=0$.
\vskip2pt
(ii) For $q=2$, assume $\gamma\in {\t C}^2 (\SV;\mathcal{B},\C)$ and denote
$\gamma_{\lambda_1,\lambda_2}(Y,Y)=a(\lambda_1-\lambda_2), \gamma_{\lambda_1,\lambda_2}(M,N)=b$. By $(d\gamma)_{\lambda_1,\lambda_2,\lambda_3}(L,Y,Y)=a(\lambda_2^2-\lambda_3^2)=0$ and $(d\gamma)_{\lambda_1,\lambda_2,\lambda_3}(L,M,N)=b(\lambda_2+\lambda_3)=0$, we deduce that $a=b=0$.
\vskip2pt

(iii) For $q=3$, assume $\gamma\in {\t C}^3 (\SV;\mathcal{B},\C)$ and denote
$\gamma_{\lambda_1,\lambda_2,\lambda_3}(Y,Y,M)=a(\lambda_1-\lambda_2), \gamma_{\lambda_1,\lambda_2}(Y,Y,N)=b(\lambda_1-\lambda_2)$. By $(d\gamma)_{\lambda_1,\lambda_2,\lambda_3,\lambda_4}(L,Y,Y,M)=
a(\lambda_2-\lambda_3)(\lambda_2+\lambda_3+\lambda_4)=0$ and $(d\gamma)_{\lambda_1,\lambda_2,\lambda_3,\lambda_4}(L,Y,Y,N)=
b(\lambda_2-\lambda_3)(\lambda_2+\lambda_3+\lambda_4)=0$, we deduce that $a=b=0$.

\vskip2pt
(iv) For $q=4$, assume $\gamma\in {\t C}^4 (\SV;\mathcal{B},\C)$ and denote
$\gamma_{\lambda_1,\lambda_2,\lambda_3,\lambda_4}(Y,Y,M,N)=a(\lambda_1-\lambda_2)$. By $(d\gamma)_{\lambda_1,\lambda_2,\lambda_3,\lambda_4,\lambda_5}(L,Y,Y,M,N)=
a(\lambda_2-\lambda_3)(\lambda_2+\lambda_3+\lambda_4+\lambda_5)=0$, we deduce that $a=0$.
\QED

%\subsection{}
By similar arguments, we have the following result.
\begin{theo} If $\mathcal{B}=\langle N\rangle, \langle N,M\rangle$, or
$\langle N,Y,M\rangle$, then
${\rm \t H}^q (\SV;\mathcal{B},\C)\cong{\rm \t H}^q (\Vir, \C)$ and ${\rm H}^q (\SV;\mathcal{B},\C)\cong {\rm H}^q (\Vir,\C)$, for $q\in\Z_+$.
\end{theo}
%\noindent{\it Proof.~}By arguments anaglous to those of \QED
\begin{coro}For $q\in\Z_+$,
${\rm \t H}^q (\SV/{\langle N,Y,M\rangle},\C)\cong{\rm \t H}^q (Vir, \C)$
and ${\rm H}^q (\SV/{\langle N,Y,M\rangle},\C)\cong {\rm H}^q (Vir,\C)$.
\end{coro}

%\subsection{}
\begin{theo} If $\mathcal{B}=\langle Y,M\rangle$, then
\begin{enumerate}[(1)]
\item ${\rm \t H}^q (\SV;\mathcal{B},\C)\cong{\rm \t H}^q (\HV, \C)$ for $q\in\Z_+$;
\item ${\rm H}^q (\SV;\mathcal{B},\C)\cong {\rm H}^q (\HV,\C)$ for $q\in\Z_+$.
 \end{enumerate}
\end{theo}
\noindent{\it Proof.~}
(1) Since $\langle Y,M\rangle$ is an ideal of $\SV$, by Remark \ref{rela}, we can deduce that ${\tilde H}^q(\SV;\langle Y,M\rangle,\C_a)={\tilde H}^q(\SV/\langle Y,M\rangle,\C_a^{\langle Y,M\rangle})$. Since $\SV/\langle Y,M\rangle\cong \mathcal{HV}$ and $\C_a^{\langle Y,M\rangle}=\C_a$, we have ${\rm \t H}^q (\SV;\langle Y,M\rangle,\C_a)\cong{\rm \t H}^q (\mathcal{HV},\C_a)$, for $q\in\Z_+$.\\
(2) Together with the similar discussion in Theorem 2.6 in \cite{YW} and (1), we can obtain the result immediately.
\QED

By the proof of Theorem \ref{th1}, we actually obtain that ${\rm \t H}^q (\SV,\C)\cong{\rm \t H}^q (\SV;\mathcal{B},\C)$ and
${\rm \t H}^q (\SV,\C)\cong{\rm H}^q (\SV;\mathcal{B},\C)$.
Then we reformulate our main result as follows.
\begin{coro}For $q\in\Z_+$, ${\rm \t H}^q (\SV,\C)\cong{\rm \t H}^q (\HV, \C)$ and ${\rm H}^q (\SV,\C)\cong {\rm H}^q (\HV,\C)$.
\end{coro}

\begin{rema}
Recall that $\SV=(\C[\partial]L\oplus\C[\partial]N)\ltimes (\C[\partial]Y\oplus\C[\partial]M)\cong \HV\ltimes (\C[\partial]Y\oplus\C[\partial]M)$. For any $\gamma\in {\t C}^q(\HV, \C)$,
define $\gamma'\in {\t C}^q(\SV, \C)$ by
$$\gamma'(X_1,\cdots,X_q)=\begin{cases}\gamma(X_1,\cdots,X_q), &X_1,\cdots,X_q\in\{L,N\};\\
0, &otherwise.
\end{cases}$$ Then the map $\gamma\mapsto \gamma'$ induces an injective homomorphism ${\t H}^q(\HV,\C)\rightarrow {\t H}^q(\SV,\C)$.
And the main theorem \ref{th1} reads that this homomorphism is actually an isomorphism.
\end{rema}

\begin{theo} If $\mathcal{B}=\langle M\rangle$, then
\begin{enumerate}[(1)]
\item ${\rm \t H}^q (\SV;\mathcal{B},\C)\cong{\rm \t H}^q (\SV, \C)$ for $q\in\Z_+$;
\item ${\rm H}^q (\SV;\mathcal{B},\C)\cong {\rm H}^q (\SV,\C)$ for $q\in\Z_+$.
 \end{enumerate}
\end{theo}
\noindent{\it Proof.~}
(1) By the Definition \ref{rcochain}, we just need to consider all solutions in Lemma \ref{le0} with $m=0$.\vskip2pt
(i) For $q=0,1$, we can obtain that ${\rm \tilde H}^0(\SV;\langle M\rangle,\C)={\rm H}^0(\SV;\langle M\rangle,\C)=\mathbb{C}$ and ${\rm \tilde H}^1(\SV;\langle M\rangle,\C)=0$ immediately by Lemma \ref{l0} and \ref{l1}.\vskip2pt
(ii) For $q=2$, we only need to consider $(k,l,m,n)=(1,0,0,1),(2,0,0,0)$ and $(0,2,0,0)$. Let $\gamma$ be an arbitrary $2$-cocycle. With the similar discussion in Lemma \ref{l5}, we can deduce that $\gamma_{\lambda_1,\lambda_2}(L, L)=\gamma_{\lambda_1,\lambda_2}(L, N)=0$ and $\gamma_{\lambda_1,\lambda_2}(Y,Y)=a(\lambda_1-\lambda_2)$, for $a\in \C$. It is not difficult to check that $(d\g)_{\lambda_1,\lambda_2,\lambda_3}(L,Y,Y)=a(\lambda_2^2-\lambda_3^2)=0,$ which implies $a=0$. Thus, we have $\t H^2(\SV;\langle M\rangle,\C)=0$.\vskip2pt
(iii) For $q=3$, similar to the proof of Lemma \ref{l6} except the case of $(k,l,m,n)=(0,2,0,1)$. We can assume that $\gamma_{\lambda_1,\lambda_2,\lambda_3}(Y,Y,N)=h(\lambda_1-\lambda_2)$. By $(d\g)_{\lambda_1,\lambda_2,\lambda_3,\lambda_4}(L,Y,Y,N)=h(\lambda_2-\lambda_3)(\lambda_2+\lambda_3+\lambda_4)=0,$ we can deduce that $h=0$, i.e., $\gamma_{\lambda_1,\lambda_2,\lambda_3}(Y,Y,N)=0$. Then we can obtain that $\t H^3(\SV;\langle M\rangle,\C_a)=\C \Phi^1\oplus\C\Phi^2\oplus\C\Phi^3$, where $\Phi^1,\Phi^2,\Phi^3$ are as shown in Lemma \ref{l6}.\vskip2pt
(iv) For $q=4$, we only need to consider $(k,l,m,n)=(2,0,0,2),(3,0,0,1),(1,2,0,1)$ and $(2,2,0,0)$. Let $\gamma$ be an arbitrary $4$-cocycle. With a similar discussion in the proof Lemma \ref{l8}, we can deduce that $\gamma_{\lambda_1,\lambda_2,\lambda_3,\lambda_4}(L,Y,Y,N)=0$ and
\begin{align*}
&\gamma_{\lambda_1,\lambda_2,\lambda_3,\lambda_4}(L,L,N,N)=(\lambda_1-\lambda_2)(\lambda_3-\lambda_4),\\
&\gamma_{\lambda_1,\lambda_2,\lambda_3,\lambda_4}(L,L,L,N)=(\lambda_1-\lambda_2)(\lambda_2-\lambda_3)(\lambda_1-\lambda_3),
\end{align*}
are cocycles but not coboundaries.

We can assume that
\begin{align*}
&\gamma_{\lambda_1,\lambda_2,\lambda_3,\lambda_4}(L,L,Y,Y)=(\lambda_1-\lambda_2)(\lambda_3-\lambda_4)(a\lambda_1+a\lambda_2+b\lambda_3+b\lambda_4),
\end{align*}
where $a,b\in\C$. Then by $$(d\gamma)_{\lambda_1,\lambda_2,\lambda_3,\lambda_4,\lambda_5}(L,L,Y,Y,N)=-2(\lambda_1-\lambda_2)(\lambda_3-\lambda_4)[a(\lambda_1+\lambda_2)+b(\lambda_3+\lambda_4+\lambda_5)]=0,$$ we have $a=b=0$, which implies $\gamma_{\lambda_1,\lambda_2,\lambda_3,\lambda_4}(L,L,Y,Y)=0.$

Thus, we can deduce that $\t H^4(\SV;\langle M\rangle,\C_a)=\C \Psi^1\oplus\C\Psi^2$, where $\Psi^1,\Psi^2$ are as shown in Lemma \ref{l8}.\vskip2pt
(v) For $q=5$, we only need to consider $(k,l,m,n)=(3,2,0,0),(1,2,0,2)$ and $(2,2,0,1)$. Let $\gamma$ be an arbitrary $5$-cocycle. We can assume that
\begin{align*}
&\gamma_{\lambda_1,\lambda_2,\lambda_3,\lambda_4, \lambda_5}(L,L,L,Y,Y)=a(\lambda_1-\lambda_2)(\lambda_1-\lambda_3)(\lambda_2-\lambda_3)(\lambda_4-\lambda_5),\\
&\gamma_{\lambda_1,\lambda_2,\lambda_3,\lambda_4, \lambda_5}(L,Y,Y,N,N)=b(\lambda_2-\lambda_3)(\lambda_4-\lambda_5),\\
&\gamma_{\lambda_1,\lambda_2,\lambda_3,\lambda_4, \lambda_5}(L,L,Y,Y,N)=(\lambda_1-\lambda_2)(\lambda_3-\lambda_4)(x\lambda_1+x\lambda_2+y\lambda_3+y\lambda_4+z\lambda_5),
\end{align*}
where $a,b,x,y,z\in \C$.

Let $\psi$ be a $4$-cochain defined by $\psi_{\lambda_1,\lambda_2,\lambda_3}(L,Y,Y,N)=-\frac{b}{2}(\lambda_2-\lambda_3)\lambda_4$. Similarly, we have$$(d\psi)_{\lambda_1,\lambda_2,\lambda_3,\lambda_4, \lambda_5}(L,Y,Y,N,N)=b(\lambda_2-\lambda_3)(\lambda_4-\lambda_5).$$ Replacing $\g$ by $\g-d\psi$, we have $\gamma_{\lambda_1,\lambda_2,\lambda_3,\lambda_4, \lambda_5}(L,Y,Y,N,N)=0$.

Let $\phi$ be a $4$-cochain defined by $\phi_{\lambda_1,\lambda_2,\lambda_3}(L,Y,Y,N)=0$ and $\phi_{\lambda_1,\lambda_2,\lambda_3}(L,L,Y,Y)=-(\lambda_1-\lambda_2)(\lambda_3-\lambda_4)(\frac{x}{2}\lambda_1+\frac{x}{2}\lambda_2+\frac{y}{2}\lambda_3+\frac{y}{2}\lambda_4)$.
Then by direct computations, we have$$(d\phi)_{\lambda_1,\lambda_2,\lambda_3,\lambda_4, \lambda_5}(L,L,Y,Y,N)=(\lambda_1-\lambda_2)(\lambda_3-\lambda_4)(x\lambda_1+x\lambda_2+y\lambda_3+y\lambda_4+y\lambda_5).$$ Replacing $\g$ by $\g-d\phi$, we have $\gamma_{\lambda_1,\lambda_2,\lambda_3,\lambda_4, \lambda_5}(L,L,Y,Y,N)=z(\lambda_1-\lambda_2)(\lambda_3-\lambda_4)\lambda_5$. Then by
\begin{align}
&(d\gamma)_{\lambda_1,\lambda_2,\lambda_3,\lambda_4, \lambda_5,\lambda_6}(L,L,Y,Y,N,N)=2z(\lambda_1-\lambda_2)(\lambda_3-\lambda_4)(\lambda_5-\lambda_6)=0,\label{LYN1}\\
&(d\gamma)_{\lambda_1,\lambda_2,\lambda_3,\lambda_4, \lambda_5,\lambda_6}(L,L,L,Y,Y,N)=2a(\lambda_1-\lambda_2)(\lambda_1-\lambda_3)(\lambda_2-\lambda_3)(\lambda_4-\lambda_5)=0,\label{LYN2}
\end{align}
we can obtain that $z=a=0$, i.e., $\gamma_{\lambda_1,\lambda_2,\lambda_3,\lambda_4, \lambda_5}(L,L,Y,Y,N)=\gamma_{\lambda_1,\lambda_2,\lambda_3,\lambda_4, \lambda_5}(L,L,L,Y,Y)=0.$
Thus, we have $\t H^5(\SV;\langle M\rangle,\C)=0$.\vskip2pt

(vi) For $q=6$, we only need to consider $(k,l,m,n)=(3,2,0,1)$ and $(2,2,0,2)$. Let $\gamma$ be an arbitrary $6$-cocycle. We can assume that
\begin{align}
&\gamma_{\lambda_1,\lambda_2,\lambda_3,\lambda_4,\lambda_5,\lambda_6}(L,L,Y,Y,N,N)=a(\lambda_1-\lambda_2)(\lambda_3-\lambda_4)(\lambda_5-\lambda_6),\label{LYN4}\\
&\gamma_{\lambda_1,\lambda_2,\lambda_3,\lambda_4,\lambda_5,\lambda_6}(L,L,L,Y,Y,N)=b(\lambda_1-\lambda_2)(\lambda_1-\lambda_3)
(\lambda_2-\lambda_3)(\lambda_4-\lambda_5),\label{LYN3}
\end{align}
for some $a,b\in\C$. By equations (\ref{LYN1}) and (\ref{LYN2}), one can easily check that (\ref{LYN4}) and (\ref{LYN3}) are coboundaries. Thus, we have $\t H^6(\SV;\langle M\rangle,\C)=0$.\vskip2pt
(vii) Obviously, we can deduce that $\t H^7(\SV;\langle M\rangle,\C)=0$ by Lemma \ref{l11}.\vskip2pt
(viii) For $q\ge8$, $\t H^q(\SV;\langle M\rangle,\C)=0$ followed by Lemma \ref{l1}.\vskip2pt
(2) By Corollary \ref{cor4.3} and Theorem \ref{th3}, we can obtain the conclusion immediately.
\QED

Note that $\SV/{\langle M\rangle}$ is a new Lie conformal algebra of rank 3, whose cohomology is also determined.
\begin{coro}For $q\in\Z_+$,
${\rm \t H}^q (\SV/{\langle M\rangle},\C)\cong{\rm \t H}^q (\SV, \C)$
and ${\rm H}^q (\SV/{\langle M\rangle},\C)\cong {\rm H}^q (\SV,\C)$.
\end{coro}
\vskip10pt

\vspace{4mm} \noindent\bf{\footnotesize Acknowledgements.}\ \rm
{\footnotesize This work was supported by the National Natural Science Foundation of China (No. 11701345, 11871421, 12171129) and the Zhejiang Provincial Natural Science Foundation of China (No. LY20A010022) and the Fundamental Research Funds for the  Central Universities (No. 22120210554).}\\
\vskip18pt \small\footnotesize
\parskip0pt\lineskip1pt

\end{document}